# OUTLIER ROBUST CORNER-PRESERVING METHODS FOR RECONSTRUCTING NOISY IMAGES[1]

By Martin Hillebrand and Christine H. Müller

*Munich University of Technology and University of Kassel*

The ability to remove a large amount of noise and the ability to preserve most structure are desirable properties of an image smoother. Unfortunately, they usually seem to be at odds with each other; one can only improve one property at the cost of the other. By combining $M$-smoothing and least-squares-trimming, the *TM*-smoother is introduced as a means to unify corner-preserving properties and outlier robustness. To identify edge- and corner-preserving properties, a new theory based on differential geometry is developed. Further, robustness concepts are transferred to image processing. In two examples, the *TM*-smoother outperforms other corner-preserving smoothers. A software package containing both the *TM*- and the *M*-smoother can be downloaded from the Internet.

**1. Introduction.** In recent years, image processing has become an important issue due to the rapid development of digitization and its applications in both industry and science. A fundamental operation in image processing is the reconstruction of a noisy digital image. A procedure denoising an image aims to achieve two objectives:

- removing as much noise as possible;
- preserving as much of the true signal as possible.

These goals are difficult to achieve at the same time and, in particular, it is difficult to remove outliers and to preserve discontinuities simultaneously.

Recently, several edge- and corner-preserving smoothing methods have been proposed. Some are methods based on wavelets and related methods (see, e.g., [2, 4, 5] and the references therein). Other methods are based on

Received January 2005; revised March 2006.

[1]Supported by the Friedrich Ebert Foundation and by Grant Mu 1031/4-1 of the Deutsche Forschungsgemeinschaft.

special local estimators where the reconstructed pixel value is calculated by pixel values in a neighborhood (window). Such a neighborhood is usually provided by a kernel function, and so these estimators are called *kernel estimators*. Chu et al. [3] proposed the use of an $M$-kernel estimator based on a redescending objective function, while Polzehl and Spokoiny [18, 19] proposed methods based on an adaptive choice of the kernel function. However, none of these methods can eliminate isolated outliers, that is, none of them is outlier robust. The methods based on wavelets and other regularization methods are even known to be peak-preserving which, in particular, means that outliers are preserved.

On the other hand, there are many reconstruction methods which are able to remove outliers. The most prominent ones in image analysis are kernel estimators based on outlier robust location estimators, such as the median smoother studied in [10] and the estimators based on least trimmed squares estimators studied by Meer et al. [11, 12], Rousseeuw and Van Aelst [23] and Müller [15, 16, 17]. But these estimators are not corner-preserving.

The ability to preserve corners and the ability to remove outliers seem to be contradictory properties. Some methods can have both properties, but not simultaneously. In these cases, one can switch from one property to the other by changing a few parameters. For regularization methods, this can be done, for example, by a high or low penalty function. For the $M$-estimator of Chu et al. [3], it depends on the scale parameter (see Figures 11 and 12).

Although a lot of literature has been published on edge- and corner-preserving smoothing, there is, surprisingly, no theoretical concept of two-dimensional discontinuities with nondifferentiable edge curves which could characterize, for example, a corner as we would identify it on the basis of a visual impression. Instead, theory has only been developed for dimension 1, where discontinuities can only be jumps, or for differentiable edge curves, as Polzehl and Spokoiny [19] have done.

This paper fills the gap—an intuitive differential geometric framework is set up in which edges and corners are properly defined. This allows for a definition of asymptotic corner- (resp. edge-) preserving as consistency at a corner (resp. edge) point.

On the other hand, there is also an absence of any formalization of the quality of removing noise in terms of robustness against irregular distributed noise such as outliers. For this purpose, we transfer robustness concepts to the image analysis context.

Having constructed such a formal framework for image smoothing, we can show that the $M$-smoother introduced by Chu et al. [3] has the remarkable property of being both asymptotically corner-preserving and robust for large samples. However, in the finite case it turns out that it is not robust against outliers.

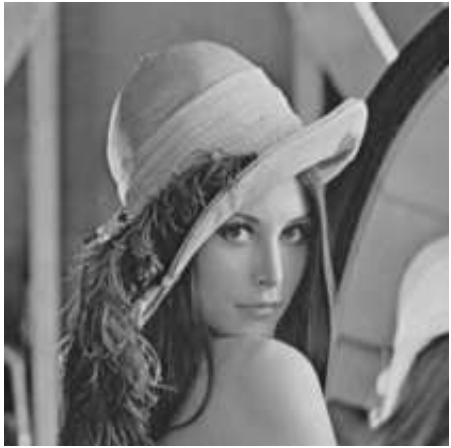

Fig. 1. *Original image.*

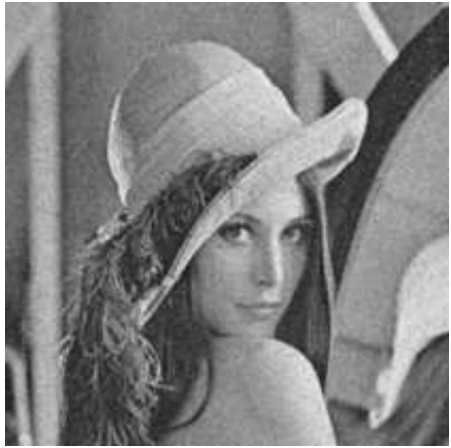

Fig. 2. *Noisy image.*

Combining the Chu smoother with a trimming procedure, we introduce the trimmed $M$-smoother ($TM$-smoother) which unifies the corner-preserving properties of the Chu smoother and an excellent outlier robustness. This is a remarkable combination of properties; we have not found such properties unified in any of the other existing methods. The $TM$-smoother is reasonably successful in distinguishing between corners and outliers. Note that the $TM$-smoother is not a two-step estimator. The $M$-smoothing part uses both the trimmed and the untrimmed data set. It is the sophisticated interplay of trimming and $M$-smoothing which gives the attractive combination of properties; the $M$-smoother alone does not remove the outliers, while the trimming procedure not only eliminates outliers but also "regular" values from corners.

The following example illustrates the outlier robustness property combined with the good smoothing property of the $TM$-smoother. To the original image (Figure 1), noise is added (Figure 2). The $TM$-smoother (Figure 3) performs better in outlier removal than the AWS-estimator of Polzehl and Spokoiny [18] (Figure 4). More details on this example are given in Section 2.

Section 2 provides more details about the $M$- and $TM$-smoothers, illustrated with further examples. Hereafter, by "$M$-smoother" we always mean the redescending $M$-smoother introduced by Chu et al. In Section 3, edges and corners are defined based on a differential geometric approach. Consistency and corner preservation are treated in Section 4. Also, model assumptions are discussed. In Section 5, asymptotic and nonasymptotic robustness concepts are transferred from location estimation to nonparametric (two-dimensional) regression. It is shown under which conditions both the $M$- and the $TM$-smoother are asymptotically robust and that the $TM$-smoother is even outlier robust in the finite case. Section 6 outlines practical aspects of the estimators and Section 7 summarizes the results.

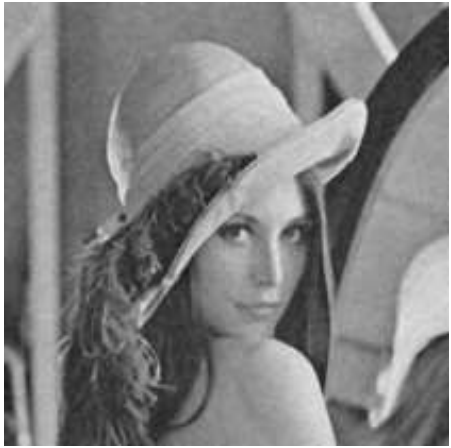
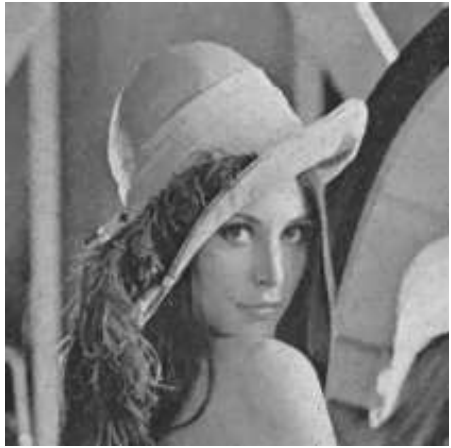

FIG. 3.   *TM-smoother.*                    FIG. 4.   *AWS.*

Since the whole theoretical framework is new, there is no standard way of showing the different properties. The proofs in the Appendix provide insights about how to work with this new theory and how it is related to standard nonparametric regression and estimation of location.

**2. *M*- and TM-estimators.**   The reason why local estimators, especially robust ones, can lose their discontinuity-preserving property when transferred from one-dimensional to two-dimensional regression can be explained as follows; see also Figures 5 and 6.

In the one-dimensional case, usually (if the jumps are not too close, which can, at least asymptotically, be assumed) the majority of the data in the neighborhood ("window") are on the "right side" of the jump. Hence, an estimator need only follow the majority of the data (as most robust estimators do) to preserve the jump. However, in the two-dimensional case, this is no longer a successful strategy—around a corner point, the majority of the data is usually on the "wrong" side of the discontinuity; see Figure 6.

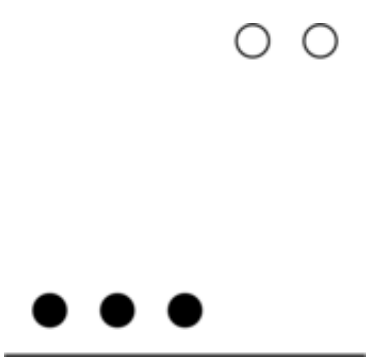
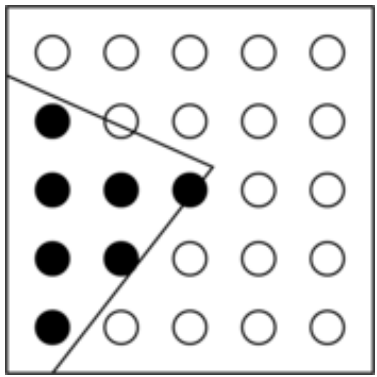

FIG. 5.   *One-dimensional discontinuity.*

FIG. 6.   *Two-dimensional discontinuity.*

The $M$-smoother of Chu et al. [3] looks for a local mode of a density estimate $H_{n,x}$ and hence also allows the estimator to be a "minority point". Therefore, it is able to preserve corners. In particular, it chooses the mode $y$ from the set of modes $\mathcal{N}_n$ of $H_{n,x}$ which is closest to the observation in the center of the window.

More formally, we consider images given by pixel values $m(x_{ij})$ (typically in a bounded interval $R$ of nonnegative numbers) at pixel positions $x_{ij}$, $i, j = 1, \ldots, n$, where we can assume without loss of generality that $x_{ij} \in [0, 1]^2$. To estimate the original image $m(x)$ on the basis of the observations $Y = (Y_{ij})_{i,j=1,\ldots,n}$, where $Y_{ij} = m(x_{ij}) + \varepsilon_{ij}$ and $\varepsilon_{ij}$ is random noise, we define the $M$-estimator of Chu et al. [3] by

$$m_n(x) := \hat{m}_{n,x}(Y) := \arg\min_{y \in \mathbb{R}} \{|y - Y_{i_0 j_0}| : y \text{ is element of } \mathcal{N}_n(x)\},$$

where

$$\mathcal{N}_n(x) := \{y \in \mathbb{R} : y \text{ is a local maximum of } H_{n,x}(y)$$
$$\text{with } y \leq Y_{i_0 j_0} \text{ if } H'_{n,x}(Y_{i_0 j_0}) \leq 0$$
$$\text{and } y > Y_{i_0 j_0} \text{ if } H'_{n,x}(Y_{i_0 j_0}) > 0\}$$

and

$$H_{n,x}(y) := \frac{1}{n^2} \sum_{i,j=1}^{n} K_{h_n}(x - x_{ij}) L_{g_n}(y - Y_{ij}).$$

Here, $(i_0, j_0) := \arg\min_{(i,j) \in \{1,\ldots,n\}^2} \|x - x_{ij}\|_2$ [if $x^k = (x_{ij}^k + x_{(i+1),j}^k)/2$ for $k = 1$ or $2$, then define $i_0 := i$ and analogously for $j_0$] and $K_{h_n}(x) := 1/h_n^2 K(x/h_n)$, $L_{g_n}(y) := 1/g_n L(y/g_n)$ with kernel functions $K : \mathbb{R}^2 \to \mathbb{R}$ and $L : \mathbb{R} \to \mathbb{R}$ and bandwidths $h_n$, $g_n \in (0, \infty)$, respectively. Since it is easier to handle zeros of a function than minima, we note that $m_n(x)$ is an element of $\{y : H'_{n,x}(y) = 0\}$. The estimator $m_n(x)$ can be calculated by means of the Newton–Raphson method starting at the center of the window $(i_0, j_0)$ and searching for the next maximum of $H_{n,x}(y)$ in the ascending direction. Existence and uniqueness of this estimator follow as in the one-dimensional case (see [8]).

One can imagine that the estimator will reach the wrong mode if the starting point $y_{i_0 j_0}$ is an outlier (see Figure 7). The basic idea of the TM-smoother is to trim the data set from which the density estimator $H_{n,x}$ is computed so that outliers cannot generate additional modes. Then the starting point remains the same, even if it is not used for the density estimation. If the starting point is an outlier, it needs to go a long way to the next mode of $H_{n,x}$, which is then hopefully the correct one (see Figure 8). We will later see that this strategy is successful.

Reducing the data set so that possible outliers are eliminated is achieved by the trimming procedure of the least-trimmed squares (LTS) estimator introduced by Rousseeuw [21] (see also [22]). Define the set of indices of observations in the window which contains all positive kernel weights by

$$J_{n,x} := \{(i,j) \in \{1, \ldots, n\}^2 : \|x - x_{ij}\|_\infty \leq h_n\}.$$

Then the *l-trimmed LTS-estimator* is defined as

$$m_{\text{LTS},l}(x) := \arg\min_{y \in \mathbb{R}} \left\{ \sum_{k=1}^{\#J_{n,x}-r} s_{(k)}(y) \right\},$$

where $(s_{(k)}(y))_{k \in \{1, \ldots, \#J_{n,x}\}}$ is the order statistic of $\{s_{ij}(y) = (y - Y_{ij})^2 : (i,j) \in J_{n,x}\}$, $l \in (0, 0.5)$ and $r := \lfloor \#J_{n,x} \cdot l \rfloor$.

Rousseeuw and Van Aelst [23] applied the LTS-estimator to image analysis, but without formalizing the two-dimensional regression model. A detailed model and a qualitative robustness analysis are provided by Müller [15, 16, 17]. However, the LTS estimator is not corner-preserving. To obtain a corner-preserving and outlier robust estimator, we do not need the LTS-estimate itself, only the trimmed set of observations

$$R_{n,l}(x) := \{(i,j) \in J_{n,x} : s_{ij}(m_{\text{LTS},l}(x)) \leq s_{(\lceil (1-l) \cdot \#J_{n,x} \rceil)}(m_{\text{LTS},l}(x))\}.$$

Then the *trimmed M-estimator* or *TM-smoother* is basically the *M*-estimator where the density estimate is based on the trimmed data set.

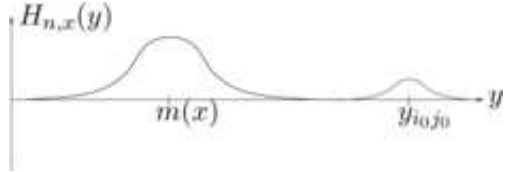

Fig. 7.   $H_{n,x}(y)$.

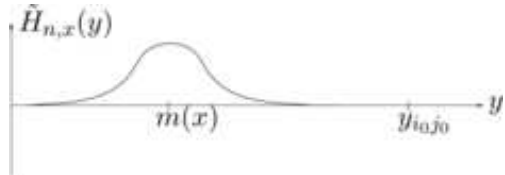

Fig. 8.   $\tilde{H}_{n,x}(y)$ based on a trimmed data set.

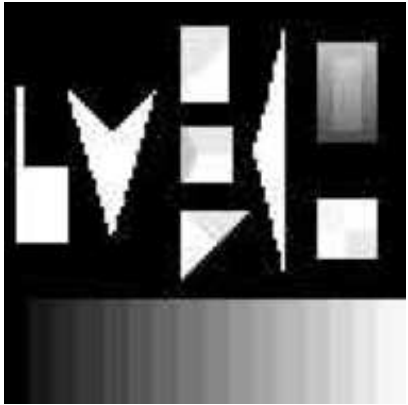

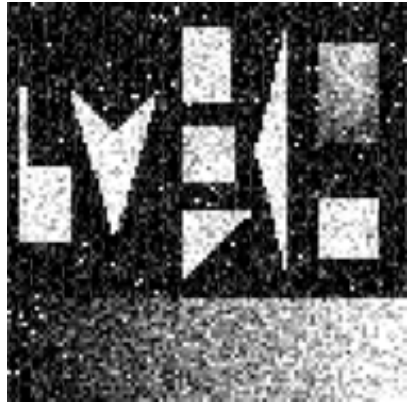

Fig. 9. *Original image.*     Fig. 10. *Noisy image.*

DEFINITION 1. The *TM-smoother* $m_{n,r}(x)$ is defined as follows:

$$m_{n,l}(x) := \hat{m}_{n,x,l}(Y)$$
$$:= \arg\min\{|y - Y_{i_0 j_0}| : y \text{ is an element of the closure of } \mathcal{N}_{n,l}(x)\},$$

where

$$\mathcal{N}_{n,l}(x) := \{y \in \mathbb{R} : y \text{ is a local maximum of } \tilde{H}_{n,x}(y) \text{ such that } \tilde{H}_{n,x}(y) > 0,$$
$$\text{with } y < Y_{i_0 j_0} \text{ if } \tilde{H}'_{n,x}(Y_{i_0 j_0}) < 0$$
$$\text{and } y > Y_{i_0 j_0} \text{ if } \tilde{H}'_{n,x}(Y_{i_0 j_0}) > 0\}$$

and

$$\tilde{H}_{n,x}(y) := \frac{1}{n^2} \sum_{(i,j) \in R_{n,l}(x)} K_{h_n}(x - x_{ij}) L_{g_n}(y - Y_{ij}).$$

$K_{h_n}$, $L_{g_n}$ and $(i_0, j_0)$ are defined as for the *M-smoother*.

How this estimator performs in practice can be seen in the following example, the "SUSAN" image given by Smith and Brady [25], downloaded from www.springerlink.com/content/?k=international+of+computer+vision. It is a $100 \times 100$ pixel image containing geometric figures with different kinds of edges and corners (see Figure 9).

To each pixel, normally distributed random background noise with a standard deviation of 26 [which is about 10% of the range of values since the brightness is linearly scaled from 0 (black) to 255 (white)] is added. In addition to the background noise (residuals) which has expectation 0 and bounded support, white colored outliers are added so that the model looks like

$$Y_{ij} = (1 - \delta_{ij})(m(x_{ij}) + \varepsilon_{ij}) + \delta_{ij} \cdot 255,$$

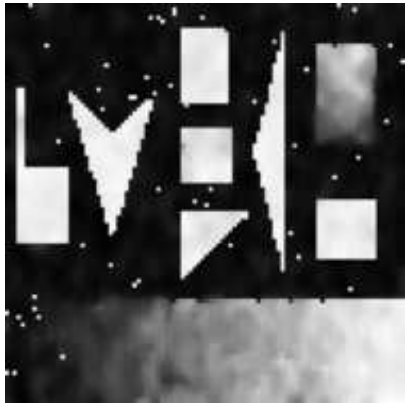
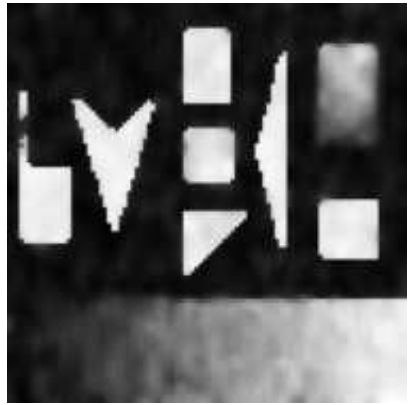

Fig. 11. *Redescending M-smoother,* $g_n = 54.5$.

Fig. 12. *Redescending M-smoother,* $g_n = 85$.

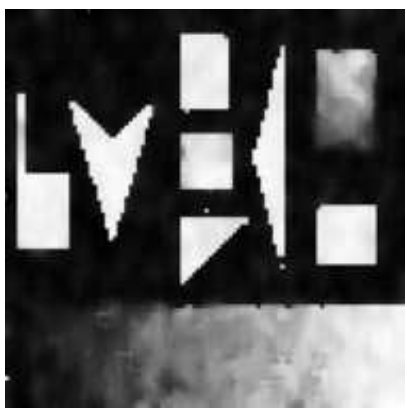
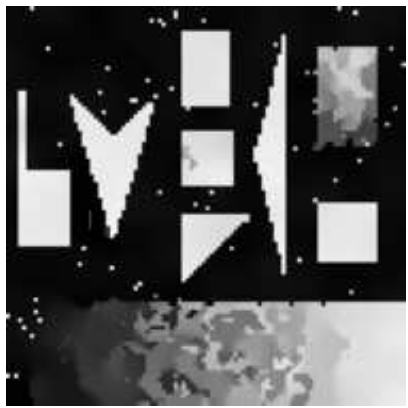

Fig. 13. *TM-smoother.*

Fig. 14. *Adaptive weights smoother.*

where $\delta_{ij}$ are i.i.d. Bernoulli distributed random variables with $p = 0.01$, in other words $\delta_{ij} \sim B(0.01)$; see Figure 10.

The noisy image is then smoothed by the *M*- or *TM*-smoother. In Figure 11, the *M*-smoother is used with parameter $g_n = 54.5$, automatically calculated as the median of the interquartile ranges. We see that corners are preserved but that outliers also are. If one increases the scale (and smoothing) parameter to $g_n = 85$, then outliers are removed, but corners are too (see Figure 12).

Applying the *TM*-smoother with $l = 0.15$ and automatically chosen scale parameter $g_n = 54.5$ to the test image in Figure 10 leads to the result in Figure 13. Now, the corners are preserved and the outliers are deleted. The software package *epsi* contains the *M*-kernel smoother and the *TM*-

TABLE 1

*Mean absolute error (MAE) and mean squared error (MSE) of the reconstructed SUSAN image*

| Method | MAE | MSE |
|---|---|---|
| Noisy image | 27.8 | 2025.0 |
| Redescending $M$-smoother, $g_n = 54.5$ | 16.1 $(-42\%)$ | 609.0 $(-70\%)$ |
| Redescending $M$-smoother, $g_n = 85$ | 19.2 $(-31\%)$ | 662.8 $(-67\%)$ |
| $TM$-smoother, $g_n = 54.5$, $l = 0.15$ | 13.9 $(-50\%)$ | 350.9 $(-83\%)$ |
| Adaptive weights smoother | 21.1 $(-24\%)$ | 795.4 $(-61\%)$ |

kernel smoother implemented in the R-library and is downloadable from cran.r-project.org.

The comparison with other corner-preserving methods shows that they are not able to delete outliers. For example, Figure 14 provides the result for the adaptive weights smoother (AWS) of Polzehl and Spokoiny [18] which appeared in their study as one of the best corner-preserving methods.

The existence of the original image gives us—in addition to the visual impression—a second criterion for the performance of an estimator: it enables us to compute the absolute and quadratic "distances" of the smoothed noisy picture from the original, the mean absolute error (MAE) $n^{-2} \sum_{i=1, j=1}^{n} |m(x_{ij}) - m_n(x_{ij})|$ and the mean squared error (MSE) $n^{-2} \sum_{i=1, j=1}^{n} (m(x_{ij}) - m_n(x_{ij}))^2$, respectively. In Table 1, the results for the different redescending $M$-kernel smoothers are given. This table also contains the corner-preserving adaptive weights smoother (AWS) of Polzehl and Spokoiny [18].

In the example given in the Introduction, we see that the $TM$-smoother also performs well on real images where the structure of the image is more complex. It is the "Lena image" which is famous in the image processing community and which can be downloaded from sipi.usc.edu/database/. To the true $512 \times 512$ pixel image, normally distributed background noise with standard deviation 17 and 1.6% outliers was added: 0.8% "salt" (white outliers) and 0.8% "pepper" (black outliers). The trimming parameter $l = 0.15$ was chosen and $h_n$ was set at 0.004, resulting in a $5 \times 5$ pixel window where three data points are trimmed. $g_n = 25.5$ was again automatically calculated. The MSE and MAE corresponding to Figures 3 and 4 are given in Table 2.

**3. Edges and corners.** While the set of discontinuities of a one-dimensional almost everywhere continuous regression function is usually the union of "jumps," the two-dimensional case is much more complicated. Here, the set of discontinuities of an a.e. continuous regression function is—apart from functions without a visual structure—a one-dimensional subset of the image that can have different shapes like the borderlines in the examples in Figure 15.

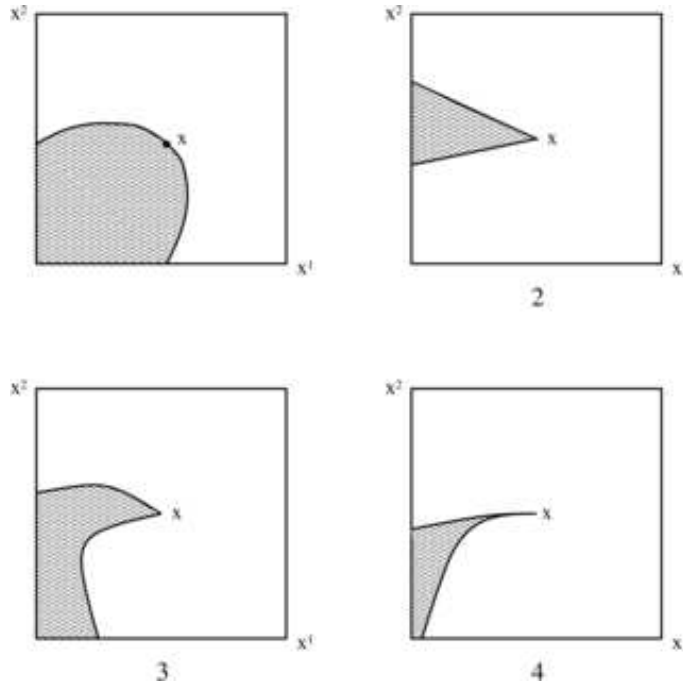

Fig. 15.    *Two-dimensional discontinuities.*

To obtain a formal characterization of the discontinuities, we turn briefly to differential geometry; see, for example, [24].

Let $I := [a, b] \subset \mathbb{R}$ be a compact interval and let $x = \binom{x^1}{x^2} : I \to \mathbb{R}^2$ be continuous. Then the set

$$\gamma := \{x(t) : t \in I\}$$

is called a *parametrically-defined (parametrized) plane curve*. The curve $\gamma$ is called *regular* if the derivatives of $x^1(t)$ and $x^2(t)$ exist. If the derivatives satisfy $\|x'(t)\| = 1$ for all $t \in I$, then the curve has a *natural* parametrization. A curve is called a *simple* or *Jordan curve with respect to given parametrization* $x = x(t)$, $t \in I$, if $x(t)$ is injective on $[a, b]$ or, if the curve is closed [i.e., $x(a) = x(b)$] on $(a, b)$.

TABLE 2
*Mean absolute error (MAE) and mean squared error (MSE) of the reconstructed Lena image*

| Method | MAE | MSE |
|---|---|---|
| Noisy image | 16.0 | 598.8 |
| *TM*-smoother, $g_n = 25.5$, $l = 0.15$ | 6.07 ($-62\%$) | 77.0 ($-87\%$) |
| Adaptive weights smoother | 6.51 ($-59\%$) | 255.9 ($-57\%$) |

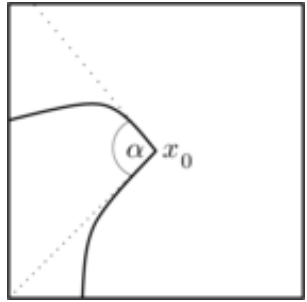

Fig. 16. *Directional tangents.*

Heretofore, we could use standard definitions. But for our highly specialized topic of interest, we have to create some special structures. For geometric singularities (points where the natural parametrization is not differentiable) we introduce the following definition.

DEFINITION 2. *If $\gamma$ is a simple curve with a natural parametrization on $I \setminus \{t_0\}$ for some $t_0 \in I$ and the limits $\lim_{t \nearrow t_0} x'(t)$ and $\lim_{t \searrow t_0} x'(t)$ exist, then the pair of directional tangents of $\gamma$ in $x_0 = x(t_0)$ is defined as*

$$T_l(\gamma, x_0) := \left\{ z \in \mathbb{R}^2 : z = x_0 - \lambda \cdot \lim_{t \nearrow t_0} x'(t), \lambda \in [0, \infty) \right\}$$

*and*

$$T_r(\gamma, x_0) := \left\{ z \in \mathbb{R}^2 : z = x_0 + \lambda \cdot \lim_{t \searrow t_0} x'(t), \lambda \in [0, \infty) \right\}.$$

Note that if $x'(t)$ is Lipschitz continuous on $I \setminus \{t_0\}$, then the directional tangents exist, by the Cauchy criterion. In Figure 16, the directional tangents, which intersect at angle $\alpha$, are represented by dotted lines.

If $x_0$ is a regular point, then the angle between the two directional tangents is $\alpha = \pi$ and the union of the directional tangents is equal to the tangent at that point. But if we have a cuspidal point (see the fourth image of Figure 15), then the directional tangents are equal and the angle between the two directional tangents is $\alpha = 0$. Hence, "real" corners in a visual sense, such as those in Images 2 and 3 of Figure 15, are characterized by the fact that the angle between the two directional tangents satisfies $\alpha \in (0, \pi) \cup (\pi, 2\pi)$.

DEFINITION 3. *Let $\gamma$ be a simple curve having parametrization*

$$x = x(t), \qquad t \in I,$$

which is natural and has a bounded second derivative $x''(t)$ in some open interval $I' \subset I$, except at a point $x_0 = x(t_0)$, $t_0 \in I'$. Then $x_0$ is called a *corner point with angle* $\alpha$ if the two directional tangents of $x_0$ intersect with angle $\alpha \in (0, \pi) \cup (\pi, 2\pi)$.

It is apparent that the corner point is well defined, that is, that the pair of directional tangents exists.

DEFINITION 4. An *edge curve* is a simple closed curve with a natural parametrization and a bounded second derivative, except at a finite number of corner points.

In the following, we will consider images for which the discontinuities can be described by edge curves; see Assumption ($\mathcal{B}2$) below. This assumption allows for the consideration of a broad variety of images and every other image can be arbitrarily well approximated by such feasible images. For the corner-preserving property, we need the rather strong condition ($\mathcal{A}2$) which requires that the discontinuities are always larger than the background noise. However, other estimates are not able to preserve corners, even if there is no noise at all. Also, we will find versions of the $M$-smoother (resp. $TM$-smoother) which are robust against a violation of the distribution assumption ($\mathcal{A}2$) [resp. ($\mathcal{A}2'$)].

The bandwidth $g_n$ is a crucial smoothing parameter—the larger $g_n$, the smoother the reconstructed image is; the smaller $g_n$, the more discontinuities are preserved. Asymptotically, Chu et al. [3] suggest $g_n \to 0$. Then $H_{n,x}$ converges to the density of the distribution of the residuals $\varepsilon_i$, as is shown in detail for the one-dimensional case in [8]. But since, in this case, even a small contamination of the residual distribution may cause a large bias in the estimator, it is not robust, as is shown in Section 5. Therefore, we suggest choosing a constant $g_n$ which attains robustness. This choice is also consistent with the automatic parameter selection in our software package. Consistency and asymptotic robustness are studied for both situations (for $g_n \to 0$ and for constant $g_n = g$). For all asymptotic results, we assume that the observations consist of the true signal and some additive noise, that is, $Y_{ij} = m(x_{ij}) + \varepsilon_{ij}$.

To prove consistency and robustness for the scale parameter converging to zero, we make the following assumptions with respect to the error distribution concerning smoothness and the number of modes ($\mathcal{A}1$). For the corner-preserving property, we need the additional assumption that the contrast at the discontinuities is larger than the noise spread ($\mathcal{A}2$).

($\mathcal{A}1$)  The regression errors (background noise) $\varepsilon_{ij}$ are independent and identically distributed with a density function $f$ supported on a bounded or unbounded interval $\mathcal{I} \subset \mathbb{R}$ such that the Lipschitz continuous derivative $f'$ has the property $f'(y) \neq 0$ for all $y \in \mathcal{I} \setminus \{0\}$ (i.e., $f$ is strongly unimodal in 0).

($\mathcal{A}2$) As Assumption ($\mathcal{A}1$), but with the additional assumption that $f$ is supported on a bounded interval $(a_1, a_2)$ such that $a_2 - a_1 < d$ (where $d$ is the jump height; see ($\mathcal{B}2$)).

Further assumptions, collectively denoted by $\mathcal{B}$, are equidistant spacing ($\mathcal{B}1$), conditions on the shape of the image ($\mathcal{B}2$), usual assumptions on the kernel function ($\mathcal{B}3$) and score function ($\mathcal{B}4$) and standard asymptotic parameter choice ($\mathcal{B}5$):

($\mathcal{B}1$) The design points are $x_{ij} := (\frac{i-1/2}{n}, \frac{j-1/2}{n})$, $i, j = 1, \ldots, n$.

($\mathcal{B}2$) The regression function is $m(x) := \mu(x) + d\mathbb{1}_D(x)$, where $m(x)$ is defined on $[0,1]^2$, $\mu(x)$ is continuous on $[0,1]^2$, $d > 0$, and $D$ is a nonempty closed set with a boundary $\partial D$ which is the disjoint union of a finite number of edge curves. Observe that a relaxation of these assumptions, allowing $d = d(x)$ to be smooth in $x$ and bounded below by some constant $d_0 > 0$, is possible.

($\mathcal{B}3$) $K(u) \geq 0$ on $(-1,1)^2$ and equals 0 elsewhere, $K(u)$ is Lipschitz continuous, $K(0) > 0$ and $\int K(u) \, du = 1$.

($\mathcal{B}4$) $L(v)$ is a nonnegative function, has a Lipschitz continuous derivative and satisfies $L(0) \neq 0$, $\int L(v) \, dv = 1$, $\int L(v)|v| \, dv < \infty$ and $\int L'(v)|v| \, dv < \infty$.

($\mathcal{B}5$) As $n \to \infty$, we have $g_n \to 0$, $h_n \to 0$, $n^{-1} h_n^{-2} \to 0$ and $n^{-1} h_n^{-1} g_n^{-2} \to 0$.

For the robust version of the estimators (fixed $g$), we can relax the assumptions ($\mathcal{A}1'$) resp. ($\mathcal{A}2'$) on $f$. The additional assumptions ($\mathcal{B}'$) differ from those of the nonrobust case only in ($\mathcal{B}4$) and ($\mathcal{B}5$) due to the fixed scale parameter.

($\mathcal{A}1'$) The regression errors $\varepsilon_{ij}$ are independently and identically distributed with density function $f$ which has bounded or unbounded support $\mathcal{I} \subset \mathbb{R}$, which is symmetric on $[-g, g]$, strictly decreasing on $(0, \infty) \cap \mathcal{I}$ and strictly increasing on $(-\infty, 0) \cap \mathcal{I}$ [i.e., $f$ is (weakly) unimodal].

($\mathcal{A}2'$) As Assumption ($\mathcal{A}1'$), but with the additional assumption that the density function $f$ is supported on the interval $(-a, a)$ and that $2a + 2g < d$.

($\mathcal{B}4'$) $L$ has two Lipschitz continuous derivatives and is nonnegative, symmetric, supported on $(-1, 1)$ and strongly unimodal on its support: $L'$ is positive on $(-1, 0)$. Finally, $L''$ has a finite number of zeros in $(-1, 1)$.

($\mathcal{B}5'$) $g_n = g$ is constant and as $n \to \infty$, $h_n \to 0$ and $n^{-1} h_n^{-1} \to 0$.

## 4. Corner preservation and consistency.

4.1. *Consistency and corner preservation for large samples.* For a location estimator, consistency is a desirable property ensuring that the estimate becomes better if the sample size increases. This concept can be transferred to nonparametric regression by calling an estimator $m_n(x)$ *consistent at* $x$ if for arbitrarily small $\varepsilon > 0$,

$$\lim_{n \to \infty} P(|m_n(x) - m(x)| > \varepsilon) = 0.$$

However, consistency depends on the shape of the regression function (in our context, the image) around $x$. Usually, an image smoother is said to be consistent if it is consistent at points with smooth neighborhoods. In the context of the class of images which we consider in our theory according to Assumption ($\mathcal{B}2$), we specify this feature as follows.

DEFINITION 5. An estimator $m_n$ is called *consistent in smooth regions* if for all $x_0 \in (0,1)^2 \setminus \partial D$ and all $\varepsilon > 0$,

$$\lim_{n \to \infty} P(|m_n(x_0) - m(x_0)| > \varepsilon) = 0.$$

For local smoothers, consistency in smooth regions can be derived from consistency of the corresponding location estimator since observations in a shrinking neighborhood are then asymptotically independently and identically distributed. Usually, the one-dimensional jump-preserving property, which is consistency at discontinuities (jumps) of a one-dimensional regression function, can be transferred to consistency at smooth edge curves, that is, edge curves without corner points (singularities), because then, at least asymptotically, the majority of the observations lies on the "correct" side of the curve. This is the case for the estimators considered in [17], for example.

Here, we show a substantially stronger consistency—the asymptotic preservation of corners. For technical reasons, we consider boundary points $x_0 \in \partial D$ which are rational, that is, $x_0 \in \partial D \cap \mathbb{Q}^2 =: \partial D_{\mathbb{Q}}$. Then there exists a subsequence $(n_l)_{l \in \mathbb{N}}$ such that $x_0$ is a grid point for all $n_l$.

DEFINITION 6. (a) An estimator $m_n$ is called *edge-preserving* for large samples if for all regular points $x_0 \in \partial D_{\mathbb{Q}}$ and for all $\varepsilon > 0$,

$$\lim_{l \to \infty} P(|m_{n_l}(x_0) - m(x_0)| > \varepsilon) = 0.$$

(b) An estimator $m_n$ is called *$\alpha$-corner-preserving* for large samples if for all corner points $x_0 \in \partial D_{\mathbb{Q}}$ with angle $\in (\alpha, \pi) \cup (\pi, 2\pi)$ and for all $\varepsilon > 0$,

$$\lim_{l \to \infty} P(|m_{n_l}(x_0) - m(x_0)| > \varepsilon) = 0.$$

(c) An estimator $m_n$ is called *corner-preserving* for large samples if it is edge-preserving and $\alpha$-corner-preserving for all $\alpha \in (0, \pi) \cup (\pi, 2\pi)$.

Therefore, only edge preservation was shown for some estimators (see, e.g., [19]). Theorem 1 shows consistency in smooth regions (a) and even, under stronger assumptions, the corner-preserving property of the $M$-smoother (b).

**THEOREM 1.** (a) *Let Assumptions ($\mathcal{A}1$) and ($\mathcal{B}$), or ($\mathcal{A}1'$) and ($\mathcal{B}'$), hold. Then the $M$-estimator of Chu et al. [3] is consistent in smooth regions.*

(b) *Let Assumptions ($\mathcal{A}2$) and ($\mathcal{B}$), or ($\mathcal{A}2'$) and ($\mathcal{B}'$), hold. Then the $M$-estimator of Chu et al. [3] is corner-preserving for large samples.*

From consistency of the $M$-smoother, uniform consistency and convergence of the integrated mean squared error can be derived as the following corollary states.

**COROLLARY 1.** *Let $m_n$ be the $M$-estimator of Chu et al. [3], let Assumptions ($\mathcal{A}1$) and ($\mathcal{B}$), or ($\mathcal{A}1'$) and ($\mathcal{B}'$), hold and let $C \subset (0,1)^2$ be a compact set and $U \subset (0,1)^2$ an open set with $\partial D \cap C \subset U$. Then for all $\varepsilon > 0$:*

(a) $\lim_{n\to\infty} P(\sup_{x\in C\setminus U} |m_n(x) - m(x)| > \varepsilon) = 0.$

*Under the additional assumptions ($\mathcal{A}2$) and ($\mathcal{A}2'$), respectively, we even have, for all $\varepsilon > 0$:*

(b) $\lim_{n\to\infty} P(\int_C |m_n(x) - m(x)|^2\, dx > \varepsilon) = 0;$

(c) $\lim_{n\to\infty} \int \int_C |m_n(x) - m(x)|^2\, dx\, dP = 0.$

Consistency of the $TM$-smoother can be derived from consistency of the $M$-smoother. However, one must take care that the right mode of the score function $\tilde{H}_{n,x}$ is not affected by the trimming procedure. The trimming proportion $l$ should not be too large. However, in real applications this is not a strong restriction since outliers are typically sparse enough that a very small $l$ can be chosen. Hence, even very small corners are preserved. To formulate the conditions for this, let $F$ be the distribution function which satisfies $F' = f$.

**THEOREM 2.** (a) *Let Assumptions ($\mathcal{A}1'$) and ($\mathcal{B}'$) hold and suppose $l < \min\{F(0), 1 - F(0)\}$. Then the $TM$-smoother $m_{n,l}$ is consistent in smooth regions.*

(b) *Let Assumptions ($\mathcal{A}2'$) and ($\mathcal{B}'$) hold and suppose $l < \frac{\alpha}{8} \cdot \min\{F(0), 1 - F(0)\}$. Then the $TM$-smoother $m_{n,l}$ is $\alpha$-corner-preserving for large samples.*

The extension of Theorem 2 to uniform convergence needs consistency of the trimming procedure, which is implied in the consistency of the LTS-estimator. To our knowledge, the consistency of the LTS-estimator has only been shown for symmetric distributions; see [1, 26]. It seems likely that it also holds for asymmetric distributions, but thus far, no one has succeeded in proving it. Hence, uniform convergence and convergence of the integrated mean squared error cannot be carried over from Corollary 1 to the *TM*-estimator yet, although we believe that it also holds.

### 4.2. *Corner preservation for finite samples.*

In the examples, we see that both the *TM*- and the *M*-smoother have good corner-preserving properties already in the finite case. This is based on the fact that the estimator is able to preserve minor features of the data sample in the window.

If we consider a corner, as sketched in Figure 6, and assume that there is no noise, that we have only two "colors" (e.g., black and white), that the pixel in the center of the window is inside the corner and that $g_n$ is sufficiently small, then such a corner is preserved by the *M*-smoother, regardless of how sharp the corner is. The *TM*-smoother requires at least $\lfloor l \cdot \#J_{n,x} \rfloor + 1$ pixels "inside" the corner to preserve it. In practice, the latter is not a strong assumption. In our examples, we used parameters allowing the preservation of all corners with more than three pixels (see Section 6 for more details).

Estimators which follow the majority of the data, like many outlier robust estimators, do not have such a strong property—in a $5 \times 5$ pixel window, they need at least 13 pixels inside the corner, which therefore must have an angle of more than $3/4\pi$.

## 5. Robustness.

### 5.1. *Large sample robustness.*

Besides the question of asymptotic "correctness" of the estimator under certain assumptions (which is answered by consistency), it is of interest to consider how the estimator is influenced by a violation of the assumptions, in particular by a contamination of the distribution of the error noise (this also includes outliers).

For estimation of location, Hampel [6] introduced large sample robustness (see also [9]). An estimator is called *robust for large samples* if a small contamination of the distribution of the observations causes only a small bias of the estimator asymptotically.

We transfer this concept from the location case to nonparametric regression. For a precise definition, we need the *Lévy metric* on the space $\mathcal{P}$ of probability measures on $\mathbb{R}$,

$$d_L(P,Q) := \min\{\varepsilon : F(y-\varepsilon) - \varepsilon \le G(y) \le F(y+\varepsilon) + \varepsilon \text{ for all } y \in \mathbb{R}\},$$

where $F$ and $G$ are the distribution functions of the probability measures $P$ and $Q$, respectively. The $\varepsilon$-*Lévy neighborhood of* $P$ is defined as

$$U_{L,\varepsilon}(P) = \{Q \in \mathcal{P} : d_L(P,Q) \le \varepsilon\}.$$

Let $m : J \subset \mathbb{R}^2 \to I \subset \mathbb{R}$, $x \mapsto m(x)$, be a regression function and let $Y := (Y_{ij})_{i,j=1,\dots,n}$, where $Y_{ij}$ are observations at $x_{ij} \in J$. For the estimator $\hat{m}_{n,x} : \mathbb{R}^{n \times n} \to \mathbb{R}$, let $(P)^{\hat{m}_{n,x}(Y)}$ be the distribution of $\hat{m}_{n,x}(Y)$ if $P$ is the distribution of the i.i.d. residuals $Y_{ij} - m(x_{ij})$.

DEFINITION 7. *The estimator $\hat{m}_{n,x}(Y)$ is called* robust for large samples *at $P$ in $x$ if for all $\varepsilon^* > 0$, there exist $\varepsilon > 0$ and $N \in \mathbb{N}$ such that*

$$d_L((P)^{\hat{m}_{n,x}(Y)}, (Q)^{\hat{m}_{n,x}(Y)}) \le \varepsilon^* \qquad \text{for all } Q \in U_{L,\varepsilon}(P) \text{ and } n \ge N.$$

Note that we use the shorter form $m_n(x)$, instead of $\hat{m}_{n,x}(Y)$, in the remainder of the paper.

The basic idea of the $M$-smoother of Chu et al. [3] is that one searches for a local mode of a density estimate $H_{n,x}$. Chu et al. suggest $g_n \to 0$ since then $H_{n,x}$ converges to the density function of the distribution of the residuals. But a small change in the noise distribution may cause an additional mode in the density function, thereby causing a fairly large bias of the estimator.

This observation led to the analysis of the asymptotics of the estimator with constant scale parameter $g_n$. In this case, $H_{n,x}$ no longer converges to the density function, rather to some other function $h$ which has the mode at the same place as the density, but which is less sensitive to contamination of the noise distribution.

The following theorems summarize these results. Their proofs give interesting insights as to how we can analyze large sample robustness properties in nonparametric regression if standard methods cannot be applied.

THEOREM 3. (a) *Let Assumptions ($\mathcal{B}$) hold and let $P$ be a distribution satisfying ($\mathcal{A}1$). Further, let $x_0 \in (0,1)^2$. Then the M-estimator $m_n(x_0)$ of Chu et al. [3] is not robust at $P$ in $x_0$ for large samples.*

(b) *Let Assumptions ($\mathcal{B}'$) hold and let $P$ be a distribution satisfying ($\mathcal{A}1'$). Let $x_0 \in (0,1)^2 \setminus \partial D$. Then the M-estimator $m_n(x_0)$ of Chu et al. [3] is pointwise robust for large samples at $P$ in $x_0$.*

If $P$ satisfies ($\mathcal{A}2'$) and $x_0$ is a grid point for some $n \in \mathbb{N}$, then the estimator is even robust at corners $x_0 \in \partial D$.

Also, according to the following theorem, the $TM$-smoother is robust for large samples.

THEOREM 4. *Let Assumptions ($\mathcal{B}'$) hold and let $P$ be a distribution satisfying ($\mathcal{A}1'$). Let $x_0 \in (0,1)^2 \setminus \partial D$. Then the TM-estimator $m_{n,l}(x_0)$ is pointwise robust for large samples at $P$ in $x_0$.*

The interesting result that large sample robustness of the $M$- and the $TM$-smoother crucially depends on the asymptotic choice of the scale parameter $g_n$ gives important information about a sensible parameter choice. The median of the interquartile ranges of the observations in the windows turns out to be a good choice for an automatic parameter selection. Moreover, it converges to some constant $g \neq 0$ and hence fits into the large sample robust modeling scheme.

5.2. *Finite sample robustness.* However, a very special kind of distribution contamination, that is, the presence of outliers, causes problems for the $M$-smoother—as can be seen in Figure 11, the $M$-estimator is not capable of removing outliers with a sensible parameter choice. Indeed, one can choose a much larger smoothing parameter $g_n$, but then corners are no longer preserved. Apparently, the asymptotic robustness property does not take effect in this case.

Hence, we need another robustness concept characterizing estimators which are able to remove outliers, such as the $TM$-smoother.

In 1971, Hampel [6] introduced a quantitative robustness measure called the *breakdown point* of an estimator. It is the minimal quota of observations which must be arbitrarily biased so that the estimator tends to $\pm\infty$. The extension of this concept to linear models can be found in Müller [14]. The special case of a breakdown point in two-dimensional nonparametric regression is treated in [16].

DEFINITION 8.   Let $x \in (h_n, 1 - h_n)^2$ and let

$$(y)_{J_{n,x}} := \{y_{ij} : (i,j) \in J_{n,x}\}$$

be the set of observations in the window $U_{h_n}(x)$. Let

$$\mathcal{Y}_{n,r,y} := \{(z)_{J_{n,x}} : z_{ij} \neq y_{ij} \text{ for at most } r \text{ of the } z_{ij}\}.$$

Then the *maximum bias of an estimator $\hat{m}_{n,x}$ by replacing $r$ observations of $(y)_{J_{n,x}}$* is defined as

$$B(\hat{m}_{n,x}, (y)_{J_{n,x}}, r) := \max\{|\hat{m}_{n,x}((y)_{J_{n,x}}) - \hat{m}_{n,x}((z)_{J_{n,x}})| : (z)_{J_{n,x}} \in \mathcal{Y}_{n,r,y}\}.$$

The *breakdown point of $\hat{m}_{n,x}$ by replacing observations of $(y)_{J_{n,x}}$* is defined as

$$\varepsilon^*(\hat{m}_{n,x}, (y)_{J_{n,x}}) := \min\left\{\frac{r}{\#J_{n,x}} : r \in \mathbb{N} \text{ with } B(\hat{m}_{n,x}, (y)_{J_{n,x}}, r) = \infty\right\}$$

and the *breakdown point of $\hat{m}_{n,x}$ by replacing observations* is defined as

$$\varepsilon^*(\hat{m}_{n,x}) := \min\{\varepsilon^*(\hat{m}_{n,x}, (y)_{J_{n,x}}) : (y)_{J_{n,x}} \in \mathbb{R}^{\#J_{n,x}}\}.$$

We can now define outlier robustness.

DEFINITION 9. An estimator $m_n(x) = \hat{m}_{n,x}$ is called *outlier robust* if its breakdown point by replacing observations is larger than $1/\#J_{n,x}$.

Although redescending $M$-estimators are known to have high breakdown points (see, e.g., [13, 20]), this does not hold for the $M$-estimator of Chu et al. [3]. This is because a high breakdown point is only achieved if the $M$-estimator is defined as the global maximum of the score function. As soon as the $M$-estimator is defined as some local maximum such as the $M$-estimator of Chu et al. [3], the high breakdown point property is lost—it is obvious that the $M$-estimator of Chu et al. [3] has a breakdown point of $1/\#J_{n,x}$. However, the $TM$-estimator is outlier robust, as we see from the following theorem.

THEOREM 5. *Let Assumptions (B3) and (B4′) hold and let $x \in (h_n, 1 - h_n)^2$.*

(a) *the $M$-estimator of Chu et al. [3] is* not *outlier robust;*

(b) *if $l \in [1/\#J_{n,x}, 1/2)$, then the $TM$-estimator $m_{n,l}$ is outlier robust and, in particular,*

$$\varepsilon^*(m_{n,l}(x)) > l.$$

This means that even if a fraction $l$ of the observations in the window are outliers, the estimator does not "break down." We believe that a stronger result also holds: that the $TM$-smoother is even consistent in the presence of a fraction of outliers not larger than $l$.

We can now clearly see the role of the parameter $l$: it removes outliers if they are at most $100 \cdot l\%$ of the observations in the window and it preserves corners if they contain more than $100 \cdot l\%$ observations; see also Section 4.2.

**6. Computational aspects.** Both the $M$-smoother and the $TM$-smoother are contained in the R package *epsi*. Since the $TM$-smoother has a considerably better robustness property and a smoothing quality which is no worse, it is preferable to use the $M$-smoother. In any case, the $M$-smoother can be regarded as a special $TM$-smoother with $l = 0$.

Like Chu et al. [3], we use the Gaussian density with mean 0 and standard deviation 1 as the kernel function $L$ and the product density of the same distribution as the kernel function $K$. Also, $L$ and $K$ are set to zero outside $[-1, 1]$ and $[-1, 1]^2$, respectively. The estimator can be used without choosing any parameters—reasonable parameters are set as default values. The window size (determined by $h_n$) turns out to be optimal at $5 \times 5$ for images of up to $10^6$ pixels which gives, for example, $h_n = 0.02$ for a $100 \times 100$ pixel image. For a $5 \times 5$ window, we use $l = 0.15$ as the trimming proportion, which means that three observations of the window are trimmed. This implies that any corner consisting of more than three pixels is preserved (cf. Figure 6). The scale smoothing parameter $g_n$ is calculated automatically as the median of the interquartile ranges. Although it is a good choice, one can sometimes improve the result slightly by adjusting the parameter manually.

The algorithm can be sketched as follows:

For each pixel,

- trim the data set in the window;
- calculate $\tilde{H}_{n,x}$ based on the trimmed data set;
- choose the starting point $Y_{i_0 j_0}$ from the complete data set and find the closest local maximum of $\tilde{H}_{n,x}$ in the ascending direction. If the gradient is zero at the starting point and $\tilde{H}_{n,x}$ is zero at this point (which is typically the case when the starting point is an outlier) then search in both directions for the closest local maximum.

For the maximization, we used a Newton–Raphson algorithm. The choice of step size is essential—we used the Amijo step size—otherwise, the convergence is poor and the algorithm is slow. The algorithm is implemented in C++. The source code can be found in the R-package *epsi*.

## 7. Conclusion.

The *TM*-smoother is a good choice for images with low level background noise, outliers, edges and corners. It is able to preserve corners and edges, and between the discontinuities, it has good smoothing properties, even if the image is not homogeneous as is desirable, for example, for the AWS-estimator of Polzehl and Spokoiny [18, 19]. However, its quality becomes worse if the variance of the background noise is too large. In that case, another method such as the AWS-estimator may be a better choice.

This article also provides a theoretical framework for image reconstruction which has not previously existed. By means of several properties, local smoothers can be compared with respect to smoothing/preservation quality and robustness. The proofs provide insight into how we can work with this theory.

## APPENDIX A: PROOFS

For the proofs of the theorems, we need the following lemmas. Particularly essential for the proof of Theorem 1 is Lemma A.1. It claims that the sum of the kernel weights of the pixel positions in $D$ converges for $x_0 \in \partial D$. For this purpose, let $U_{h_n}(x_0) := \{x \in [0,1]^2 : \|x_0 - x\|_\infty \leq h_n\}$ be the window around $x_0$ with respect to $h_n$ and let $\bar{G}_n(x_0) := D \cap U_{h_n}(x_0)$. If $x_0 \in (0,1) \setminus \partial D$, then $\bar{G}_n(x_0) = \varnothing$ or $\bar{G}_n(x_0) = U_{h_n}(x_0)$ for sufficiently large $n$. If $x_0 \in \partial D$, then $\varnothing \neq \bar{G}_n(x_0) \subsetneq U_{h_n}(x_0)$ for all $n \in \mathbb{N}$ (see Figure 17).

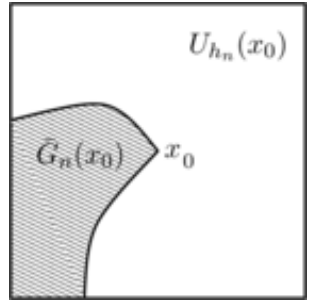

FIG. 17.    $\bar{G}_n(x_0) \subset U_{h_n}(x_0)$.

LEMMA A.1.    *Let $x_0 \in \partial D$ and let Assumptions (B1), (B2), (B3) and (B5), or Assumptions (B1), (B2), (B3) and (B5′), hold. Then there exists $G(x_0) \subset [-1,1]^2$ such that*

$$\frac{1}{n^2} \sum_{x_{ij} \in \bar{G}_n(x_0)} K_{h_n}(x_0 - x_{ij}) = \int_{G(x_0)} K(u)\,du + o(1)$$

*and $1 > \int_{G(x_0)} K(u)\,du > 0$.*

PROOF.    We prove Lemma A.1 only for corner points $x_0$. It is apparent that the proof for consistency at corner points also holds for regular points. For some fixed $n_0 \in \mathbb{N}$, let the set of discontinuities $\partial D \cap U_{h_{n_0}}(x_0)$ be described by the edge curve $x(t)$ and let $t_0 \in I$ be such that $x(t_0) = x_0$. In the following proof, we always assume $n \ge n_0$ and $n_0$ is sufficiently large so that $x_0$ is the only corner point in $\partial D \cap U_{h_{n_0}}$.

Let $b_l := \lim_{t \nearrow t_0} x'(t)$ and $b_r := \lim_{t \searrow t_0} x'(t)$. Then

$$T_l(\gamma, x_0) = \{z \in \mathbb{R}^2 : z = x_0 - \lambda \cdot b_l, \lambda \in [0, \infty)\}$$

and

$$T_r(\gamma, x_0) = \{z \in \mathbb{R}^2 : z = x_0 + \lambda \cdot b_r, \lambda \in [0, \infty)\}.$$

Consider the rotation of parameters $\Theta : \mathbb{R}^2 \to \mathbb{R}^2$ defined by

$$x \mapsto \tilde{x} := \begin{pmatrix} c^1 & c^2 \\ -c^2 & c^1 \end{pmatrix} x,$$

where

$$c = \frac{b_l + b_r}{\|b_l + b_r\|_2}.$$

Recall that $\|b_l\|_2 = \|b_r\|_2 = 1$ because $x(t)$ is a natural parametrization. $\Theta$ maps $c$ (which is the normalized sum of the direction vectors $b_l, b_r$ of the directional tangents of $x_0$ onto the $\tilde{x}^1$-axis $\binom{1}{0}$; see Figure 18.

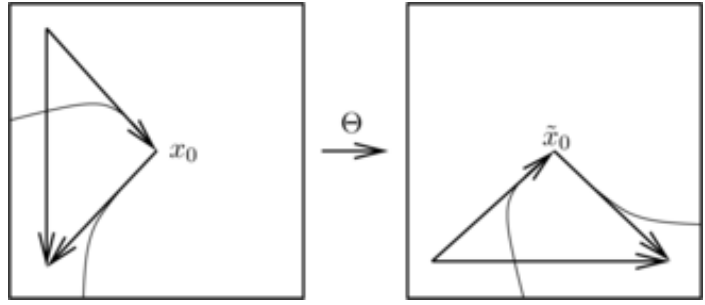

Fig. 18.   *Rotation $\Theta$.*

Observe that $\tilde{b}_l^1 > 0$ and $\tilde{b}_r^1 > 0$. This can be seen as follows. By the Cauchy–Schwarz inequality,

$$\langle b_r, b_r + b_l \rangle = \|b_r\|_2^2 + \langle b_r, b_l \rangle \geq \|b_r\|_2^2 - |\langle b_r, b_l \rangle| > 0.$$

Hence,

$$\langle b_r, c \rangle = \left\langle b_r, \frac{1}{\|b_r + b_l\|_2}(b_r + b_l) \right\rangle > 0$$

and since $\Theta$ is a rotation,

$$\tilde{b}_r^1 = \left\langle \tilde{b}_r, \begin{pmatrix} 1 \\ 0 \end{pmatrix} \right\rangle = \langle \Theta(b_r), \Theta(c) \rangle = \langle b_r, c \rangle > 0.$$

$\tilde{b}_l^1 > 0$ is shown analogously.

This, together with the Lipschitz continuity of $x'(t)$, implies that there exists a neighborhood $U_\varepsilon(t_0)$ of $t_0$ such that $\tilde{x}^{1\prime}(t) > 0$ on $U_\varepsilon(t_0)$ and hence $\tilde{x}^1$ is invertible. Then there exists, with $\tilde{U}_1 := (\tilde{x}^1)^{-1}(U_\varepsilon(t_0))$, a function $g : \tilde{U}_1 \to \mathbb{R}$ such that $g(\tilde{x}^1(t)) = \tilde{x}^2(t)$ for all $t \in U_\varepsilon(t_0)$ and which is twice differentiable on $\tilde{U}_1 \setminus \{\tilde{x}_0^1\}$.

The function $g$ can be given explicitly as

$$g(z) = \tilde{x}^2((\tilde{x}^1)^{-1}(z))$$

for $z \in \tilde{U}_1$. Hence, for $z \in \tilde{U}_1 \setminus \{\tilde{x}_0^1\}$,

$$g'(z) = \frac{(\tilde{x}^2)'((\tilde{x}^1)^{-1}(z))}{(\tilde{x}^1)'((\tilde{x}^1)^{-1}(z))}$$

and

$$\lim_{z \nearrow \tilde{x}_0^1} g'(z) = \frac{\tilde{b}_l^2}{\tilde{b}_l^1} =: \beta_l, \qquad \lim_{z \searrow \tilde{x}_0^1} g'(z) = \frac{\tilde{b}_r^2}{\tilde{b}_r^1} =: \beta_r.$$

Since the curve is simple, there exists, for sufficiently small $\tilde{U}_2$, a neighborhood $\tilde{U}_2 \subset \mathbb{R}$ such that

$$\{\tilde{x}(t) : t \in I\} \cap (\tilde{U}_1 \times \tilde{U}_2) = \{(\tilde{x}^1, g(\tilde{x}^1)) : \tilde{x}^1 \in \tilde{U}_1\}$$

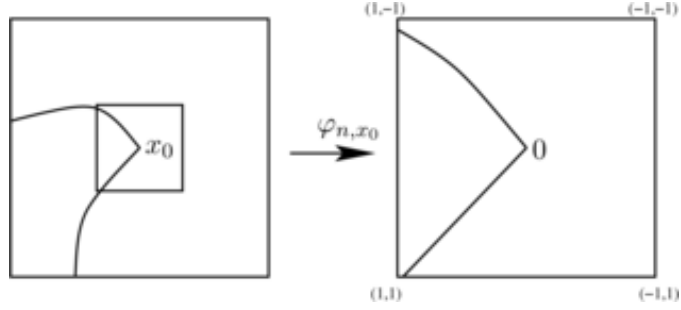

Fig. 19. $\varphi_{n,x_0}$.

and $\tilde{x}_0^2$ lies in the interior of $\tilde{U}_2$. Without loss of generality, assume that $\Theta(D) \cap (\tilde{U}_1 \times \tilde{U}_2)$ lies beneath $g$, that is, $\Theta(D) \cap (\tilde{U}_1 \times \tilde{U}_2) = \{\tilde{x} \in \tilde{U}_1 \times \tilde{U}_2 : \tilde{x}^2 \leq g(\tilde{x}^1)\}$. Then there exists $n_1 \geq n_0$ such that $\Theta(U_{h_n}(x_0)) \subset \tilde{U}_1 \times \tilde{U}_2$ for all $n \geq n_1$ and hence $\bar{G}_n(x_0) = D \cap U_{h_n}(x_0) = \{u \in U_{h_n}(x_0) : \tilde{u}^2 \leq g(\tilde{u}^1)\}$.

Moreover, there exist two Taylor expansions of $g$ at $\tilde{x}_0$,

$$\tilde{x}^2 = g(\tilde{x}^1) = \tilde{x}_0^2 + (\tilde{x}^1 - \tilde{x}_0^1)\beta_l + (\tilde{x}^1 - \tilde{x}_0^1)\eta_l(\tilde{x}^1 - \tilde{x}_0^1) \qquad \text{for } \tilde{x}^1 \leq \tilde{x}_0^1$$

and

$$\tilde{x}^2 = g(\tilde{x}^1) = \tilde{x}_0^2 + (\tilde{x}^1 - \tilde{x}_0^1)\beta_r + (\tilde{x}^1 - \tilde{x}_0^1)\eta_r(\tilde{x}^1 - \tilde{x}_0^1) \qquad \text{for } \tilde{x}^1 \geq \tilde{x}_0^1,$$

where

$$\lim_{a \to 0} \eta_i(a) = 0 \qquad \text{for } i = l, r.$$

Define the transformation

$$\varphi_{n,x_0} : U_{h_n}(x_0) \longrightarrow [-1, 1]^2$$

$$u \longmapsto \frac{1}{h_n}(x_0 - u).$$

$\varphi$ maps the window $U_{h_n}(x_0)$, which contains the support of the kernel function, onto the (mirror) unit square; see Figure 19.

Now, define

$$\bar{B}_n(x_0) := \{u \in U_{h_n}(x_0) : \tilde{u}^2 \leq \tilde{x}_0^2 + (\tilde{u}^1 - \tilde{x}_0^1)(\beta_l \mathbb{1}_{(-\infty, \tilde{x}_0^1]}(\tilde{u}^1) + \beta_r \mathbb{1}_{(\tilde{x}_0^1, \infty)}(\tilde{u}^1))\}.$$

$\bar{B}_n(x_0)$ is the area which lies, with respect to the rotated axes, "beneath" the directional tangents of $x_0$; see Figure 20.

Further, we have

$$\varphi_{n,x_0}(\bar{B}_n(x_0)) = \{u \in [-1, 1]^2 : x_0 - h_n u \in \bar{B}_n(x_0)\}$$

$$= \{u \in [-1, 1]^2 : \tilde{u}^2 \geq \tilde{u}^1 \cdot [\beta_l \mathbb{1}_{[0, \infty)}(\tilde{u}^1) + \beta_r \mathbb{1}_{(-\infty, 0)}(\tilde{u}^1)]\}.$$

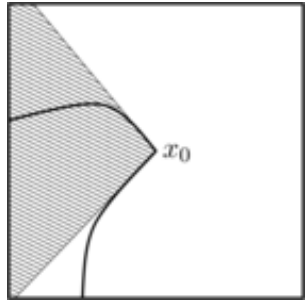

Fig. 20. $\bar{B}_n(x_0)$.

Since $\varphi_{n,x_0}(\bar{B}_n(x_0))$ is independent of $n$, we can rename it as

$$G(x_0) := \varphi_{n,x_0}(\bar{B}_n(x_0));$$

see Figure 21. Now, consider, with the Taylor expansions mentioned above,

$$G_n(x_0) := \varphi_{n,x_0}(\bar{G}_n(x_0)) = \{u \in [-1,1]^2 : \Theta(x_0 - h_n u)^2 \leq g(\Theta(x_0 - h_n u)^1)\}$$

$$= \{u \in [-1,1]^2 : \tilde{u}^2 \geq \tilde{u}^1 \cdot [(\beta_l + \eta_l(-h_n\tilde{u}^1))\mathbb{1}_{[0,\infty)}(\tilde{u}^1)$$

$$+ (\beta_r + \eta_r(-h_n\tilde{u}^1))\mathbb{1}_{(-\infty,0)}(\tilde{u}^1)]\}.$$

Define

$$\eta_{\max,n} := \Big\{ \max_{u \in [-2h_n, 2h_n]} |\eta_l(u)|, \max_{u \in [-2h_n, 2h_n]} |\eta_r(u)| \Big\}.$$

Since

$$G_n(x_0) \bigtriangleup G(x_0)$$

$$\subset \{u \in [-1,1]^2 : |\tilde{u}^2 - \tilde{u}^1(\beta_l\mathbb{1}_{[0,\infty)}(\tilde{u}^1) + \beta_r\mathbb{1}_{(-\infty,0)}(\tilde{u}^1))| \leq \eta_{\max,n}\},$$

where the *symmetric difference* is defined as $A \bigtriangleup B := (A \setminus B) \cup (B \setminus A)$, the Lebesgue measure of the symmetric difference can be estimated by $\lambda(G_n(x_0) \bigtriangleup G(x_0)) \leq 6\eta_{\max,n} = o(1)$ as $n \to \infty$. It follows immediately that

$$\int_{G_n(x_0)} K(u)\, du = \int_{G(x_0)} K(u)\, du + o(1)$$

since $K$ is bounded. Hence, it suffices to show that

$$\lim_{n \to \infty} \Big( \frac{1}{n^2} \sum_{x_{ij} \in G_n(x_0)} K_{h_n}(x_0 - x_{ij}) - \int_{G_n(x_0)} K(u)\, du \Big) = 0.$$

For the proof of this property and other technical details of the proof, see [7]. □

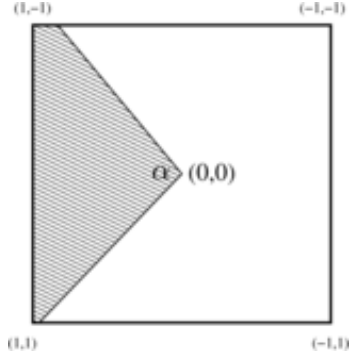

Fig. 21. $G(x_0)$.

We define the set of indices in $J_{n,x_0} = \{(i,j) \in \{1, \ldots, n\}^2 : \|x_0 - x_{ij}\|_\infty \leq h_n\}$ corresponding to $D$ by

$$I_n^{\bar{G}_n}(x_0) := \{(i,j) \in \{1, \ldots, n\}^2 : x_{ij} \in \bar{G}_n(x_0)\}.$$

Observe that for all $(i,j) \in J_{n,x_0} \setminus I_n^{\bar{G}_n}(x_0)$, we have $m(x_{ij}) = \mu(x_{ij})$ and for all $(i,j) \in I_n^{\bar{G}_n}(x_0)$, we have $m(x_{ij}) = \mu(x_{ij}) + d$.

We then have the following corollary.

COROLLARY A.1.

$$\frac{1}{n^2} \sum_{(i,j) \in J_{n,x_0} \setminus I_n^{\bar{G}_n}(x_0)} K_{h_n}(x_0 - x_{ij}) = 1 - \int_{G(x_0)} K(u)\, du + o(1).$$

Note that the equalities in Lemma A.1 and Corollary A.1 also hold for $x_0 \in (0,1)^2 \setminus \partial D$. If $x_0 \in D \setminus \partial D$, then $\int_{G(x_0)} K(u)\, du = 1$ and if $x_0 \in (0,1)^2 \setminus D$, then $\int_{G(x_0)} K(u)\, du = 0$.

Define

$$\nu_{x_0} := \int_{G(x_0)} K(u)\, du$$

and for the case that the scale parameter is converging to zero,

$$f_{d,\nu_{x_0}}(y) := \begin{cases} \nu_{x_0} f(y) + (1 - \nu_{x_0}) f(y + d), & \text{for } \nu_{x_0} \in (0,1), \\ f(y), & \text{for } \nu_{x_0} = 1 \text{ or } \nu_{x_0} = 0. \end{cases}$$

For the case that the scale parameter is fixed by $g_n = 1$, define

$$h_{d,\nu_{x_0}}(y) := \begin{cases} \nu_{x_0} h(y) + (1 - \nu_{x_0}) h(y + d), & \text{for } \nu_{x_0} \in (0,1), \\ h(y), & \text{for } \nu_{x_0} = 1 \text{ or } \nu_{x_0} = 0, \end{cases}$$

where $h(y) := \int L(y - u) P(du)$.

LEMMA A.2.   *Let $x_0 \in (0,1)^2$.*

(a) *Let Assumptions ($A1$) and ($B$) hold. Then*

$$\sup_{y \in \mathbb{R}} |EH'_{n,x_0}(y) - f'_{d,\nu_{x_0}}(y - m(x_0))| = o(1).$$

(b) *Let Assumptions ($A1'$) and ($B'$) hold. Then*

$$\sup_{y \in \mathbb{R}} |EH'_{n,x_0}(y) - h'_{d,\nu_{x_0}}(y - m(x_0))| = o(1).$$

PROOF.   (a) We provide the proof only for the case $x_0 \in \partial D$. The proof for $x_0 \in (0,1)^2 \setminus \partial D$ is the same, but even more simple. From Lemma A.1, Corollary A.1 and the Lipschitz continuity of $f'$, we obtain

$$\sup_{y \in \mathbb{R}} \left| \frac{1}{n^2} \sum_{i,j=1}^n K_{h_n}(x_0 - x_{ij}) E \frac{d}{dy} L_{g_n}(y - Y_{ij}) - f'_{d,\nu_{x_0}}(y - m(x_0)) \right|$$

$$= \sup_{y \in \mathbb{R}} \left| \frac{1}{n^2} \sum_{(i,j) \in I_n^{\bar{G}_n}(x_0)} K_{h_n}(x_0 - x_{ij}) \int \frac{d}{dy} \frac{1}{g_n} L\left( \frac{y - m(x_{ij}) - u}{g_n} \right) f(u) \, du \right.$$

$$- \int_{G(x_0)} K(u) \, du \, f'(y - m(x_0))$$

$$+ \frac{1}{n^2} \sum_{(i,j) \in J_{n,x_0} \setminus I_n^{\bar{G}_n}(x_0)} K_{h_n}(x_0 - x_{ij})$$

$$\times \int \frac{d}{dy} \frac{1}{g_n} L\left( \frac{y - \mu(x_{ij}) - u}{g_n} \right) f(u) \, du$$

$$\text{(A.1)} \hspace{4cm} \left. - \left( 1 - \int_{G(x_0)} K(u) \, du \right) f'(y - \mu(x_0)) \right|$$

$$\leq \sup_{y \in \mathbb{R}} \left\{ \frac{1}{n^2} \sum_{(i,j) \in I_n^{\bar{G}_n}(x_0)} K_{h_n}(x_0 - x_{ij}) \right.$$

$$\times \int L(v) |f'(y - m(x_{ij}) - vg_n) - f'(y - m(x_0))| \, dv$$

$$+ \frac{1}{n^2} \sum_{(i,j) \in J_{n,x_0} \setminus I_n^{\bar{G}_n}(x_0)} K_{h_n}(x_0 - x_{ij})$$

$$\times \int L(v) |f'(y - \mu(x_{ij}) - vg_n) - f'(y - \mu(x_0))| \, dv \Bigg\} + o(1)$$

$$= o(1).$$

The proof of assertion (b) is based on an equality which is analogous to equality (A.1) in the proof of (a). Note that we now have $E\frac{d}{dy}L(y - Y_{ij}) = h'(y - m(x_{ij}))$ for $(i,j) \in I_n^{\bar{G}_n}(x_0)$ and $E\frac{d}{dy}L(y - Y_{ij}) = h'(y - \mu(x_{ij}))$ for $(i,j) \in J_{n,x_0} \setminus I_n^{\bar{G}_n}(x_0)$.  $\square$

LEMMA A.3.  *Let $x_0 \in (0,1)^2$ and let Assumptions $(\mathcal{A}1)$ and $(\mathcal{B})$, or Assumptions $(\mathcal{A}1')$ and $(\mathcal{B}')$, hold. Then*

$$\lim_{n \to \infty} P\Big(\sup_{y \in \mathbb{R}} |H'_{n,x_0}(y) - EH'_{n,x_0}(y)| < \varepsilon\Big) = 1 \qquad \text{for all } \varepsilon > 0.$$

PROOF.  The proof is analogous to that of the one-dimensional case given by Hillebrand and Müller [8] in Lemma 4 using the Fourier transform of $L'$. The only difference is that $\varphi_n(u)$ must be defined as $\varphi_n(u) = n^{-2}h_n^{-2}\sum_{k,j=1}^n K(\frac{x - x_{kj}}{h_n})e^{-iuY_{kj}}$ instead of $\varphi_n(u) = n^{-1}h_n^{-1}\sum_{k=1}^n \cdot K(\frac{x - x_k}{h_n})e^{-iuY_k}$, where $i = \sqrt{-1}$. Then the condition $n^{-1}h_n^{-1}g_n^{-2} \to 0$ of Assumption $(\mathcal{B}5)$ is used instead of $n^{-1}h_n^{-1}g_n^{-4} \to 0$ for $g_n$ converging to zero. If $g_n$ is fixed, then it is clear that we only need Assumption $(\mathcal{B}5')$. Then the result can also be shown without the Fourier transform: since $L'$ is bounded by Assumption $(\mathcal{B}4')$, we obtain pointwise convergence by using Chebyshev's inequality and Corollary A.1. The Lipschitz continuity of $L'$ and $h'$ then imply the uniform convergence.  $\square$

PROOF OF THEOREM 1.  The proof of Theorem 1(a) is analogous to that of the one-dimensional case given by Hillebrand and Müller [8], replacing $\lambda$ by $\nu_{x_0}$. In particular, it is based on Lemma A.2(a) and Lemma A.3. For fixed $g_n$, the proof is the same as for $g_n \to 0$ if $f$ is replaced by $h$, Lemma A.2(b) is used instead of Lemma A.2(a) and $f_{d,\nu_{x_0}}$ is replaced by $h_{d,\nu_{x_0}}$. We see that $h_{d,\nu_{x_0}}$ has the same properties as $f_{d,\nu_{x_0}}$ because of Assumptions $(\mathcal{A}2')$ and $(\mathcal{B}4')$. For this purpose, note, in particular, that the support of $h$ is $(-a - g, a + g)$ and that $h$ is strongly unimodal.  $\square$

From the proofs of Lemma A.2 and Lemma A.3, we obtain Lemma A.4.

LEMMA A.4.  *Let $C \subset (0,1)^2$ be a compact set and $U \subset (0,1)^2$ an open set with $\partial D \cap C \subset U$.*

(a) *Under Assumptions $(\mathcal{A}1)$ and $(\mathcal{B})$, we have*

$$\lim_{n \to \infty} P\Big(\sup_{x \in C \setminus U} \sup_{y \in \mathbb{R}} |H'_{n,x}(y) - f(y - m(x))| < \varepsilon\Big) = 1 \qquad \text{for all } \varepsilon > 0.$$

(b) *Under Assumptions $(\mathcal{A}1')$ and $(\mathcal{B}')$, we have*

$$\lim_{n \to \infty} P\Big(\sup_{x \in C \setminus U} \sup_{y \in \mathbb{R}} |H'_{n,x}(y) - h(y - m(x))| < \varepsilon\Big) = 1 \qquad \text{for all } \varepsilon > 0.$$

PROOF OF COROLLARY 1.    Assertion (a) follows from the uniform convergence of $H_{n,x}$ to $f(y-m(x))$ and $h(y-m(x))$, respectively, both of which are strongly unimodal, where $h(y) := \int L(y-u)P(du)$ (see Lemma A.4). Assertions (b) and (c) hold since the Lebesgue measure of $\partial D$ is zero, so an open set $U \supset \partial D \cap C$ can be found with arbitrarily small Lebesgue measure. Then under Assumptions $(\mathcal{A}2)$ and $(\mathcal{A}2')$, respectively, because of the bounded support of the error, $|m_n(x) - m(x)|$ is also bounded, so $\int_U |m_n(x) - m(x)|^2 \, dx$ can be chosen arbitrarily small by means of an appropriate $U$. This also holds if $m_n(x)$ does not converge to $m(x)$ on $U$. □

LEMMA A.5.    Let be $y_{u,n} = \min\{Y_{ij}; (i,j) \in R_{n,l}(x)\}$ and $y_{o,n} = \max\{Y_{ij}; (i,j) \in R_{n,l}(x)\}$. Then

   (a)  $\tilde{H}'_{n,x}(y) = H'_{n,x}(y)$ for $y \in (y_{u,n} + g, y_{o,n} - g)$ if $y_{u,n} + g \leq y_{o,n} - g$;

   (b)  $\tilde{H}'_{n,x}(y) \geq H'_{n,x}(y)$ for $y \in [y_{u,n}, y_{u,n} + g]$ and $\tilde{H}'_{n,x}(y) \leq H'_{n,x}(y)$ for $y \in [y_{o,n} - g, y_{o,n}]$;

   (c)  $\tilde{H}'_{n,x}(y) > 0$ for $y \in (y_{u,n} - g, y_{u,n})$ and $\tilde{H}'_{n,x}(y) < 0$ for $y \in (y_{o,n}, y_{o,n} + g)$;

   (d)  $\tilde{H}_{n,x}(y) = 0$ and $\tilde{H}'_{n,x}(y) = 0$ for $y \in \mathbb{R} \setminus (y_{u,n} - g, y_{o,n} + g)$.

PROOF OF THEOREM 2.    The proof of Theorem 1 is based on the fact (cf. [8]) that for all $\varepsilon_1, \varepsilon' > 0$, there exist $C_1$, $C_2$ and $n_0 \in \mathbb{N}$ such that

$$P(C_1 \leq Y_{i_0,j_0} - m(x_0) \leq C_2) \geq 1 - \varepsilon_1 \qquad \text{for } n \geq n_0,$$

$$H'_{n,x}(y) > 0 \qquad \text{on } [m(x_0) + C_1, m(x_0) - \varepsilon'],$$

$$H'_{n,x}(y) < 0 \qquad \text{on } [m(x_0) + \varepsilon', m(x_0) + C_2].$$

Let $q_{l,n}$ be the $l$-quantile and $q_{1-l,n}$ the $(1-l)$-quantile of $\{Y_{i,j}; (i,j) \in J_{n,x_0}\}$. Since quantiles are asymptotically linear, we have (see [17])

$$\lim_{n \to \infty} q_{l,n} = q_l \quad \text{and} \quad \lim_{n \to \infty} q_{1-l,n} = q_{1-l},$$

where $q_l$ and $q_{1-l}$ are the quantiles of the distribution given by $f_{d,\nu_{x_0}}(y - m(x_0))$. Since

$$\lim_{n \to \infty} \frac{\# I_n^{\tilde{G}_n}(x_0)}{\# J_{n,x}} = \frac{\lambda(G(x_0))}{\lambda([-1,1]^2)} = \nu_{x_0} \geq \frac{\alpha}{8} \qquad \text{in the case } x_0 \in \partial D,$$

there exists for every $x_0 \in (0,1)^2$ and for $\varepsilon'$ sufficiently small, some $n_1 \geq n_0$ such that for all $n \geq n_1$,

$$y_{u,n} \leq q_{l,n} < m(x_0) - \varepsilon' < m(x_0) + \varepsilon' \leq q_{1-l,n} \leq y_{o,n}.$$

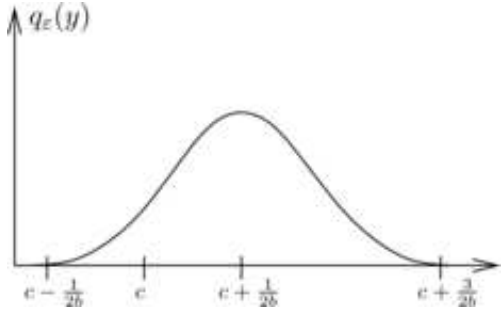

Fig. 22.   $q_\varepsilon(y)$.

By Lemma A.5, we have for all $n \geq n_1$,

$$\tilde{H}'_{n,x}(y) > 0 \qquad \text{on } [\max\{y_{u,n}, m(x_0) + C_1\}, \ m(x_0) - \varepsilon'],$$

$$\tilde{H}'_{n,x}(y) < 0 \qquad \text{on } [m(x_0) + \varepsilon', \ \min\{y_{o,n}, m(x_0) + C_2\}].$$

Hence, because $Y_{i_0,j_0} \in [m(x_0) + C_1, m(x_0) + C_2]$ with probability greater than $1 - \varepsilon_1$, the closest local maximum of $\tilde{H}_{n,x}$ to $Y_{i_0,j_0}$ with $\tilde{H}_{n,x}(y) \neq 0$ lies in $[m(x_0) - \varepsilon', m(x_0) + \varepsilon']$.   $\square$

PROOF OF THEOREM 3(a).   It suffices to show the claim for $x \in (0,1)^2 \setminus \partial D$. For arbitrarily small $\varepsilon > 0$, we will create a distribution which lies in the $\varepsilon$-Lévy-neighborhood of $P$ and has a multimodal density.

Let $c > 0$ be such that $\int_c^\infty f(y)\, dy > 0$. Further, let $\delta := -f'(c) > 0$. Consider

$$q_\varepsilon(y) := \begin{cases} a\left(1 - \left(y - c - \dfrac{1}{2b}\right)^2 b^2\right)^2, & \text{if } y \in \left[c - \dfrac{1}{2b}, c + \dfrac{3}{2b}\right], \\ 0, & \text{otherwise}, \end{cases}$$

where $a := \sqrt{\dfrac{5\delta}{8\varepsilon}}$ and $b := \sqrt{\dfrac{32\delta}{45\varepsilon}}$; see Figure 22.

It is easily verified that $q_\varepsilon$ is continuously differentiable and Lipschitz continuous, satisfying $q'_\varepsilon(c) = \dfrac{\delta}{\varepsilon}$ and $\int q_\varepsilon(u)\, du = 1$. Hence,

$$f_\varepsilon(y) := (1 - \varepsilon) f(y) + \varepsilon q_\varepsilon(y)$$

is a density function with $f'_\varepsilon(c) = \varepsilon \cdot \delta > 0$ and the corresponding distribution $P_\varepsilon$ lies in the $\varepsilon$-Lévy-neighborhood of $P$ since

$$|F(y) - F_\varepsilon(y)| = \varepsilon \cdot |F(y) - G_\varepsilon(y)| \leq \varepsilon,$$

where $G_\varepsilon(y)$ is the distribution function of the distribution $Q_\varepsilon$ with density $q_\varepsilon(y)$ and $F_\varepsilon(y)$ is the distribution function of the distribution $P_\varepsilon$.

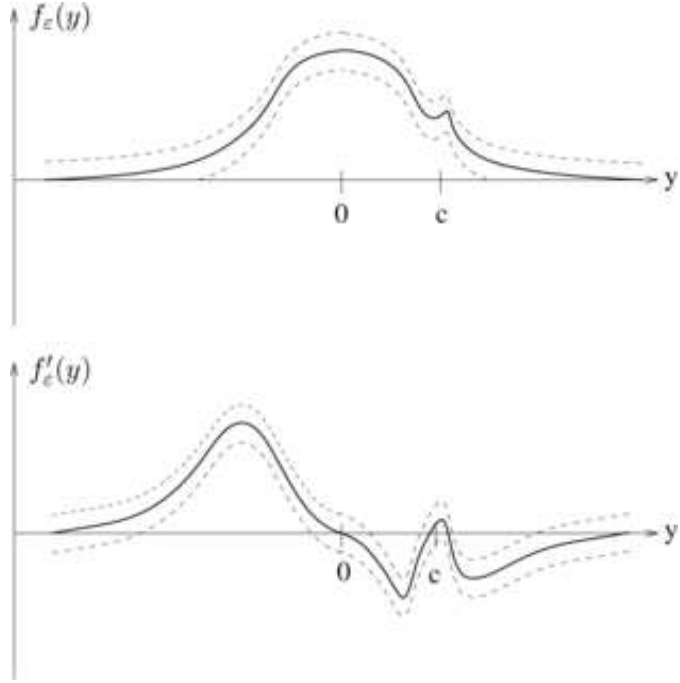

Fig. 23.   $f_\varepsilon(y)$ and $f'_\varepsilon(y)$.

Note that $f'_\varepsilon(c + 3/(2b)) < 0$ since $q'_\varepsilon(c + 3/(2b)) = 0$. Since $f_\varepsilon$ is differentiable, it has a local maximum between $c$ and $c + 3/(2b)$; see Figure 23. For sufficiently small $\varepsilon > 0$, $c + 3/(2b)$ is close to $c$ and hence

$$\int_{c + \frac{3}{2b}}^\infty f(u)\, du > 0.$$

Since Lemma A.2(a) and Lemma A.3 also hold for $f_\varepsilon(y)$, $H_{n,x}(y)$ has a local maximum in $[m(x) + c, m(x) + 3/2b)]$ with probability tending to one as $n \to \infty$. If, additionally, the starting point is greater than $m(x) + c + 3/(2b)$, then $m_n(x)$ will be greater than $m(x) + c$.

Let $(Q_\varepsilon)^{m_n(x)}$ denote the distribution of the estimator $m_n(x)$ if $Q_\varepsilon$ is the distribution of the residuals. Then if $\varepsilon_1 \geq 0$ is the probability (vanishing as $n \to \infty$) that $H_{n,x}(y)$ has no local maximum in $[m(x) + c, m(x) + c + 3/(2b)]$, we have

$$(Q_\varepsilon)^{m_n(x)}([m(x) + c, \infty]) \geq \int_{c + \frac{3}{2b}}^\infty f_\varepsilon(u)\, du - \varepsilon_1.$$

Since by Theorem 1 we also have

$$(P)^{m_n(x)}([m(x) + c/2, \infty]) \leq \varepsilon_2$$

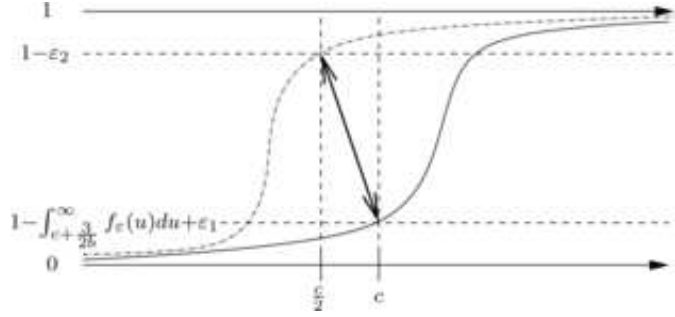

Fig. 24. *Distribution functions of $(P)^{m_n(x)}$ and $(Q_\varepsilon)^{m_n(x)}$.*

for some $\varepsilon_2 > 0$ which vanishes as $n$ becomes large, we have, as shown in Figure 24,

$$d_L((P)^{m_n(x)}, (Q_\varepsilon)^{m_n(x)}) \geq \min\left\{\int_{c+\frac{3}{2b}}^\infty f_\varepsilon(u)\,du - \varepsilon_1 - \varepsilon_2,\ \frac{c}{2}\right\}$$

$$\geq \min\left\{\int_{c+\frac{3}{2b}}^\infty f(u)\,du - \varepsilon - \varepsilon_1 - \varepsilon_2,\ \frac{c}{2}\right\}. \qquad \square$$

PROOF OF THEOREM 3(b). Let $Q_\varepsilon \in U_{L,\varepsilon}(P)$ and let $G_\varepsilon$ be its distribution function. Further, let $f_{\max} := \max_{y \in \mathbb{R}} f(y)$ and $h'_{G_\varepsilon}(y) = \int L'(y-u)\,dG_\varepsilon(u)$. Because

$$F(y) - f_{\max} \cdot \varepsilon - \varepsilon \leq F(y-\varepsilon) - \varepsilon \leq G_\varepsilon(y) \leq F(y+\varepsilon) + \varepsilon \leq F(y) + f_{\max} \cdot \varepsilon + \varepsilon,$$

we have

(A.2) $$|G_\varepsilon(y) - F(y)| \leq f_{\max} \cdot \varepsilon + \varepsilon$$

for all $y \in \mathbb{R}$. Assumption ($\mathcal{B}4'$) then implies

$$|h'_{G_\varepsilon}(y) - h'(y)| = \left|\int_{-g}^g L''(u)(G_\varepsilon(y-u) - F(y-u))\,du\right|$$

$$\leq \int_{-g}^g |L''(u)||G_\varepsilon(y-u) - F(y-u)|\,du$$

$$\leq \int_{-g}^g |L''(u)|(f_{\max} \cdot \varepsilon + \varepsilon)\,du$$

$$= C \cdot \varepsilon,$$

where $C := \int_{-g}^g |L''(u)|\,du(f_{\max} + 1)$.

Let $\varepsilon_1 > 0$ be arbitrarily small. Let $\delta := \min\{|h'(y)| : y \in [-a, -\varepsilon_1] \cup [\varepsilon_1, a]\}$. Obviously, $\delta > 0$.

Let $\varepsilon < \frac{1}{C} \cdot \frac{\delta}{2}$. Then

$$\sup_{y \in \mathbb{R}} |h'(y) - h'_{G_\varepsilon}(y)| < \frac{\delta}{2}.$$

Since Lemma A.2(b) and Lemma A.3 also hold for $G_\varepsilon$, we obtain that for arbitrarily small $\varepsilon_2 > 0$, there exists $n_0 \in \mathbb{N}$ such that with probability $1 - \varepsilon_2$,

$$\sup_{y \in \mathbb{R}} |H'_{n,x}(y) - h'_{G_\varepsilon}(y - m(x))| < \frac{\delta}{2}$$

for all $n \geq n_0$. Hence, with probability $1 - \varepsilon_2$,

$$\sup_{y \in \mathbb{R}} |H'_{n,x}(y) - h'(y - m(x))| < \delta$$

for all $n \geq n_0$. This implies:

(1) $H'_{n,x}(y) > 0$ on $[m(x) - a, m(x) - \varepsilon_1]$ and
    $H'_{n,x}(y) < 0$ on $[m(x) + \varepsilon_1, m(x) + a]$;
(2) at least one zero of $H'_{n,x}(y)$, which is a local minimum of $-H_{n,x}(y)$, lies in the $\varepsilon_1$-neighborhood of $m(x)$.

We conclude that if the starting point lies in $(m(x_{i_0}) - a, m(x_{i_0}) + a)$, the closest zero of $H'_{n,x}(y)$ in the direction searched lies, for $n \geq n_0$, in $[m(x) - \varepsilon_1, m(x) + \varepsilon_1]$, with probability larger than $1 - \varepsilon_2$. From (A.2), we have that the probability of the starting point lying in $(m(x_{i_0}) - a, m(x_{i_0}) + a)$ is greater than $1 - 2(f_{\max} + 1)\varepsilon$. Hence,

$$(Q_\varepsilon)^{m_n(x)}([m(x) - \varepsilon_1, m(x) + \varepsilon_1]) \geq 1 - \varepsilon_2 - 2(f_{\max} + 1)\varepsilon;$$

see Figure 25. Since by Theorem 1(a), we have

$$(P)^{m_n(x)}([m(x) - \varepsilon_1, m(x) + \varepsilon_1]) \geq 1 - \varepsilon_2,$$

it follows that for $n \geq n_0$,

$$d_L((P)^{m_n(x)}, (Q_\varepsilon)^{m_n(x)}) \leq \max\{2\varepsilon_1, \varepsilon_2 + 2(f_{\max} + 1)\varepsilon\}. \qquad \square$$

PROOF OF THEOREM 4.   The assertion can be seen to follow from Theorem 3(b) with arguments similar to those used to show that Theorem 2 follows from Theorem 1.   $\square$

PROOF OF THEOREM 5(b).   Let $(y)_{J_{n,x}} \in \mathbb{R}^{\#J_{n,x}}$ and set

$$y_{\min} := \min\{y_{ij} : (i,j) \in J_{n,x}\}$$

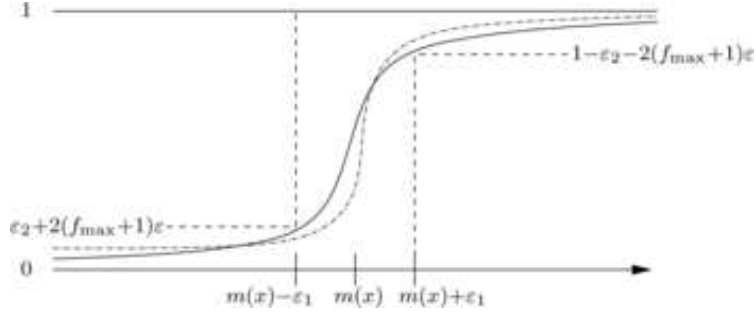

Fig. 25. *Distribution functions of $(P)^{m_n(x)}$ and $(Q_\varepsilon)^{m_n(x)}$.*

and

$$y_{\max} := \max\{y_{ij} : (i,j) \in J_{n,x}\}.$$

Let $(z)_{J_{n,x}} \in \mathcal{Y}_{n,r,y}$. Since at least $\#J_{n,x} - r$ elements of $(z)_{J_{n,x}}$ are contained in $[y_{\min}, y_{\max}]$, we have

$$(A.3) \qquad \min_{y \in \mathbb{R}} \sum_{k=1}^{\#J_{n,x}-r} s_{(k)}(y) \leq (\#J_{n,x} - r)(y_{\max} - y_{\min})^2.$$

Let $\hat{y} \in \arg\min_{y \in \mathbb{R}} \sum_{k=1}^{\#J_{n,x}-r} s_{(k)}(y)$. Then

$$\hat{y} \in \left[ y_{\min} - \sqrt{\#J_{n,x} - r}(y_{\max} - y_{\min}), y_{\max} + \sqrt{\#J_{n,x} - r}(y_{\max} - y_{\min}) \right]$$

since otherwise there is at least one $z_{i_0 j_0}$ with $z_{i_0 j_0} = y_{i_0 j_0} \in [y_{\min}, y_{\max}]$ and

$$s_{i_0 j_0}(\hat{y}) = (y_{i_0 j_0} - \hat{y})^2 > (\#J_{n,x} - r)(y_{\max} - y_{\min})^2,$$

which contradicts (A.3). If some

$$z_{i_1 j_1} \in \mathbb{R} \setminus \left[ y_{\min} - 2\sqrt{\#J_{n,x} - r}(y_{\max} - y_{\min}), y_{\max} + 2\sqrt{\#J_{n,x} - r}(y_{\max} - y_{\min}) \right],$$

then $s_{i_1 j_1}(\hat{y}) = (z_{i_1 j_1} - \hat{y})^2 > (\#J_{n,x} - r)(y_{\max} - y_{\min})^2$ and hence $(i_1, j_1) \notin R_{n,l}(x)$.

This means that all $z_{ij}$ with $(i,j) \in R_{n,l}(x)$ lie in

$$\left[ y_{\min} - 2\sqrt{\#J_{n,x} - r}(y_{\max} - y_{\min}), y_{\max} + 2\sqrt{\#J_{n,x} - r}(y_{\max} - y_{\min}) \right].$$

From the definition of $m_{n,l}(x)$, it immediately follows that $m_{n,l}(x) = \hat{m}_{n,x,l}(z)$ lies in the support of $\tilde{H}_{n,x}(z)$, which is not greater than

$$\left[ y_{\min} - 2\sqrt{\#J_{n,x} - r}(y_{\max} - y_{\min}) - g, y_{\max} + 2\sqrt{\#J_{n,x} - r}(y_{\max} - y_{\min}) + g \right].$$

This proves the claim. $\square$

**Acknowledgments.** We would like to thank Professor Chi Kang Chu for his immediate responses and explanations and the referees for their helpful comments which considerably improved the paper. We are also very grateful for Tim Garlipp's help in creating the R-package.

## REFERENCES

[1] BEDNARSKI, T. and CLARKE, B. R. (1993). Trimmed likelihood estimation of location and scale of the normal distribution. *Austral. J. Statist.* **35** 141–153. MR1244844

[2] CANDÈS, E. J. and DONOHO, D. L. (2000). Ridgelets: A key to higher-dimensional intermittency. In *Wavelets: The Key to Intermittent Information?* (B. Silvermann and J. Vassilicos, eds.) 111–127. Oxford Univ. Press.

[3] CHU, C. K., GLAD, I. K., GODTLIEBSEN, F. and MARRON, J. S. (1998). Edge-preserving smoothers for image processing (with discussion). *J. Amer. Statist. Assoc.* **93** 526–556. MR1631321

[4] DONOHO, D. L. (1999). Wedgelets: Nearly minimax estimation of edges. *Ann. Statist.* **27** 859–897. MR1724034

[5] DONOHO, D. L., JOHNSTONE, I. M., KERKYACHARIAN, G., and PICARD, D. (1995). Wavelet shrinkage: Asymptopia? (with discussion). *J. Roy. Statist. Soc. Ser. B* **57** 301–369. MR1323344

[6] HAMPEL, F. R. (1971). A general qualitative definition of robustness. *Ann. Math. Statist.* **42** 1887–1896. MR0301858.

[7] HILLEBRAND, M. (2003). On robust corner-preserving smoothing in image processing. Ph.D. dissertation, Univ. Oldenburg, Germany. Available at docserver.bis.uni-oldenburg.de/publikationen/dissertation/2003/hilonr03/hilonr03.html.

[8] HILLEBRAND, M. and MÜLLER, CH. H. (2006). On consistency of redescending $M$-kernel smoothers. *Metrika* **63** 71–90. MR2230421

[9] HUBER, P. (1981). *Robust Statistics.* Wiley, New York. MR0606374

[10] KOCH, I. (1996). On the asymptotic performance of median smoothers in image analysis and nonparametric regression. *Ann. Statist.* **24** 1648–1666. MR1416654

[11] MEER, P., MINTZ, D. and ROSENFELD, A. (1990). Least median of squares based robust analysis of image structure. In *Proc. DARPA Image Understanding Workshop* 231–254. Morgan Kaufmann, San Francisco.

[12] MEER, P., MINTZ, D., ROSENFELD, A. and KIM, D.Y. (1991). Robust regression methods for computer vision: A review. *Internat. J. Computer Vision* **6** 59–70.

[13] MIZERA, I. and MÜLLER, CH. H. (1999). Breakdown points and variation exponents of robust $M$-estimators in linear models. *Ann. Statist.* **27** 1164–1177. MR1740116

[14] MÜLLER, CH. H. (1997). *Robust Planning and Analysis of Experiments. Lecture Notes in Statistics* **124**. Springer, New York. MR1454843

[15] MÜLLER, CH. H. (1999). On the use of high breakdown point estimators in the image analysis. *Tatra Mt. Math. Publ.* **17** 283–293. MR1737714

[16] MÜLLER, CH. H. (2002). Comparison of high-breakdown-point estimators for image denoising. *Allg. Stat. Arch.* **86** 307–321. MR1924458

[17] MÜLLER, CH. H. (2002). Robust estimators for estimating discontinuous functions. *Metrika* **55** 99–109. MR1903286

[18] POLZEHL, J. and SPOKOINY, V. G. (2000). Adaptive weights smoothing with applications to image restoration. *J. R. Stat. Soc. Ser. B Stat. Methodol.* **62** 335–354. MR1749543

[19] POLZEHL, J. and SPOKOINY, V. (2003). Image denoising: Pointwise adaptive approach. *Ann. Statist.* **31** 30–57. MR1962499

[20] Riedel, M. (1984). Comparison of break points of estimators. In *Robustness of Statistical Methods and Nonparametric Statistics* (D. Rasch and M. L. Tiku, eds.) 113–116. Reidel, Dordrecht. MR0857842

[21] Rousseeuw, P. J. (1984). Least median of squares regression. *J. Amer. Statist. Assoc.* **79** 871–880. MR0770281

[22] Rousseeuw, P. J. and Leroy, A. M. (1987). *Robust Regression and Outlier Detection*. Wiley, New York. MR0914792

[23] Rousseeuw, P. J. and Van Aelst, S. (1999). Positive-breakdown robust methods in computer vision. In *Computing Science and Statistics. Models, Predictions and Computing. Proc. 31st Symposium on the Interface* (K. Berk and M. Pourahmadi, eds.) 451–460. Interface Foundation of North America, Fairfax Station, VA.

[24] Shikin, E. V. (1995). *Handbook and Atlas of Curves.* CRC Press, Boca Raton, FL. MR1412799

[25] Smith, S. M. and Brady, J. M. (1995). SUSAN—A new approach to low level image processing. *International J. Computer Vision* **23** 45–78.

[26] Yohai, V. J. and Maronna, R. A. (1976). Location estimators based on linear combinations of modified order statistics. *Comm. Statist. Theory Methods* **5** 481–486. MR0436459

Center for Mathematical Sciences
Munich University of Technology
Boltzmannstrasse 3
85747 Garching
Germany
E-mail: mhi@ma.tum.de

Department of Mathematics
and Computer Science
University of Kassel
Heinrich-Plett-Strasse 40
D-34132 Kassel
Germany
E-mail: cmueller@mathematik.uni-kassel.de